\documentclass[11pt,a4paper]{article}

\usepackage[centertags]{amsmath}
\usepackage{amsfonts}
\usepackage{amssymb}
\usepackage{epsfig}
\usepackage{hyperref}
\usepackage{color}
\usepackage{ntheorem}
\usepackage{subfigure}

\theoremstyle{break} \newtheorem{theorem}{Theorem}[section]
\theoremstyle{nonumberbreak} 
\theoremstyle{break} \newtheorem{definition}[theorem]{Definition}       
\theoremstyle{break} \newtheorem{lemma}[theorem]{Lemma}
\theoremstyle{break} \newtheorem{corollary}[theorem]{Corollary}
\theoremstyle{break} \newtheorem{example}[theorem]{Example}
\theoremstyle{nonumberbreak} \newtheorem{example1}{Example}
\theoremstyle{break} 
\theoremstyle{break} \newtheorem{remarks}[theorem]{Remarks}
\theoremstyle{break} \newtheorem{remark}[theorem]{Remark}
\theoremstyle{nonumberbreak} \newtheorem{remark1}{Remark}
\theoremstyle{nonumberbreak} \newtheorem{warning}{Warning}
\theoremstyle{break} 
\theoremstyle{break} 
\theoremstyle{break} 
\theoremstyle{break} \newtheorem{proposition}[theorem]{Proposition}

\numberwithin{equation}{section}

\newcommand{\N}{{\mathbb{N}}}

\newcommand{\R}{{\mathbb{R}}}
\newcommand{\D}{{\mathbb{D}}}
\newcommand{\C}{{\mathbb{C}}}

\renewcommand{\H}{{\mathbb{H}}}

\def\d{\mathop{{\rm d}_{\lambda}}}

\def\Re{\mathop{{\rm Re}}}
\def\Im{\mathop{{\rm Im}}}

\newcommand{\eps}{\varepsilon}

\begin{document}

\setcounter{page}{1}

\vspace*{0.0cm}
\begin{center}
 { \LARGE  \bf
 Conformal Metrics  \\[2mm]}
  {\large   Daniela Kraus and Oliver Roth}\\[2mm]
{\small University of W\"urzburg\\
Department of Mathematics \\[0.5mm]     
 D--97074 W\"urzburg\\ Germany}\\[1mm]
\today
 \end{center}

\setcounter{footnote}{1}
\footnotetext{D.K.~was supported by a HWP scholarship and 
O.R.~received partial support  from the German--Israeli Foundation (grant G--809--234.6/2003). Both authors were also supported
by a DFG grant (RO 3462/3--1).} 
 \vspace{3mm}

\section{Introduction} 
Conformal metrics connect complex analysis, differential geometry and partial
differential equations. They were used  by Schwarz \cite{Sch1891},
Poincar\'e \cite{Poi1898}, Picard \cite{Pic1890,Pic1893,Pic1905}
and Bieberbach \cite{Bie12,Bie16}, but  it was recognized by
 Ahlfors \cite{Ahl38} and Heins  \cite{Hei62} that they are {\it ubiquitous}
in complex analysis and geometric function theory. They  have been
instrumental in  important recent results such as the determination of the
spherical Bloch constant by Bonk and Eremenko \cite{BE}. 

\smallskip

The present notes grew out of a series of lectures given at the CMFT Workshop 2008 in Guwahati, India. The goal of the lectures was
to give a self--contained introduction to the rich field of conformal metrics 
assuming only some basic knowledge in complex analysis.

\smallskip

Our point of departure is the fundamental concept of {\it conformal
  invariance}. This  leads in a natural way to the notion of curvature
in Section \ref{sec:curvature} and in particular to a consideration of
  metrics with constant curvature  and
  their governing equation, the Liouville equation. The study of this equation
  and its numerous ramifications in complex analysis occupies large parts of
  our discussion. We focus on the case of negative curvature and
first find in an elementary and constructive way all {\it radially
  symmetric} solutions and thereby explicit formulas for the hyperbolic metric of all circularly
  symmetric domains like disks, annuli and punctured disks.

\smallskip

In Section 3, we discuss  Ahlfors' Lemma \cite{Ahl38} including
a simplified proof of the case of equality using only Green's Theorem
and a number of  standard applications such as Pick's Theorem and
  Liouville's Theorem.  Following  a suggestion of A.~Beardon and D.~Minda
 Ahlfors' Lemma is called the {\it Fundamental Theorem} throughout this
 paper -- a grandiose title but, as we shall see, one that is fully justified.
In a next step, 
we describe the
elementary and elegant construction of a conformal metric on the twice--punctured plane
 with curvature bounded above by a negative constant due to Minda and Schober
  \cite{MinScho3}, which is based on earlier work of Robinson \cite{Rob39}. The power of  this tool is illustrated by proving the little Picard 
Theorem, Huber's  Theorem and Picard's Big Theorem.

\smallskip

The topic of Section 4 is M.~Heins' \cite{Hei62} celebrated theory of {\it
  SK--metrics}, which we develop for simplicity  only for 
{\it continuous}  metrics. 
Following Heins' point of view, we emphasize the analogy between SK--metrics
  and  
  subharmonic functions and prove for instance  a Gluing Lemma for SK--metrics.
 A crucial step is the solution of the Dirichlet problem
for conformal metrics of constant negative curvature  on {\it disks}. As there
  is no analog of the Poisson formula, the construction is a little more
  involved and is carried out in detail in the Appendix using only some basic facts
  from classical potential theory and a fixed point argument. After this
  preparation it is quite easy to carry over a number of fundamental properties
and concepts from the theory of harmonic and subharmonic functions such as Harnack's
  monotone convergence theorem and Perron families to SK--metrics. As an
  application the existence of the hyperbolic metric on any domain with at
  least two boundary points is established. 
The usefulness of the Gluing Lemma is illustrated by deriving
 a precise asymptotic estimate for the hyperbolic metric close to
  an isolated boundary point. This then implies the {\it completeness} of 
the  hyperbolic metric and leads to a simple proof of Montel's Big Theorem.

\smallskip

In Section 5, we describe the general (local) solution of Liouville's equation in terms of
bounded analytic functions. This exhibits an intimate relation between 
 constantly curved conformal metrics and analytic function theory by
 associating to every metric with curvature $-1$ its {\it developing map}, which in general is a multivalued  analytic function. The connection is realized via the 
{\it Schwarzian derivative} of a conformal
metric, which makes it possible to reconstruct a constantly curved metric from its Schwarzian using the Schwarzian
differential equation. We give two applications to the hyperbolic
metric. First, a quick proof that the Schwarzian of the
hyperbolic metric of any domain with an isolated boundary point
has always a pole of order $2$  at this specific point is obtained.
Second, it is  shown that the inverse of the developing map of the
hyperbolic metric is always single--valued. This
combined with the results of Section 4  immediately  proves the Uniformization Theorem for planar domains. 
We close the paper by
showing how the techniques of the present paper can be  used to  
derive an explicit formula for the hyperbolic metric on the twice--punctured
plane originally due to Agard. This involves some simple facts from special
function theory, in particular about hypergeometric functions.

\smallskip

This paper contains no new results at all. Perhaps some of the
proofs might be considered as novel, which is more or less  a byproduct of our attempts to try to
present the material in a way as simple as possible. In order not to interrupt
the presentation, we have not included references in the main text. 
Some references can be found in notes at the end of every section, which also
contain some suggestions for further reading. We apologize for any omission.
We have included a number of exercises and we intend to provide  solutions for them at {\it www.mathematik.uni-wuerzburg.de/$\sim$kraus}.
There is no attempt to present
a comprehensive treatment of ``conformal metrics'' in this paper; we rather focus on some selected
topics and strongly recommend the monograph \cite{KL} of L.~Keen and N.~Lakic
and the 
forthcoming book of A.~Beardon and D.~Minda on the
same topic for extensive discussions of conformal metrics. This is required reading for anybody interested
in conformal metrics! The adventurous reader is directed to Heins' paper \cite{Hei62}.

\vspace{-2mm}
\subsection*{Acknowledgments}

\vspace*{-3mm}
We would like to thank all the participants of the CMFT Workshop for their
valuable questions and comments. We owe a particular debt of gratitude to
the other ``Resource persons'' of the workshop
 Lisa Lorentzen, Frode R{\o}nning, Richard Fournier, S.~Ponnusamy and
Rasa and J\"orn Steuding. We are very grateful to Alan Beardon, Edward Crane, David
Minda, Eric Schippers and Toshi Sugawa for many helpful discussions and conversations, and in
particular for being tireless advocates of the theory of conformal metrics.
And finally, we wish to say a big thank--you to Meenaxi Bhattacharjee and
Stephan Ruscheweyh.

\medskip

\begin{quote} 
{\sl \large ``Pardon me for writing such a long letter; I had not the time to write a
  short one.''}

\smallskip

\hfill{\large Lord Chesterfield}
\end{quote}

\newpage

{\bf \Large Glossary of Notation}

\vspace{1cm}

\begin{tabular}{llr}
{\bf Symbol} \hspace*{1cm} & {\bf Meaning} & {\bf Page} \\[2mm]
$\N$          & set of non--negative integers & \pageref{def:N}\\[1mm]
$\R$          & set of real numbers & \pageref{ref:R}\\[1mm]
$\C$          & complex plane & \pageref{def:G}\\[1mm]
$\C'$         & punctured complex plane $\C \backslash\{ 0 \}$ & \pageref{def:C'}\\[1mm]
$\C''$        & twice--punctured plane $\C \backslash \{0,1\}$ & \pageref{def:C''}\\[1mm]
$D$, $G$      & domains in the complex plane  & \pageref{def:G} \\[1mm]
$\overline{M}$ & closure of $M \subseteq \C$ relative to $\C$& \pageref{def:closure}\\[1mm]
$\partial M$ &  boundary of $M \subseteq \C$ relative to $\C$ & \pageref{def:p}\\[1mm]
$K_R(z_0)$    & open disk in $\C$, center $z_0 \in \C$, radius $R>0$ & 
\pageref{def:dik}\\[1mm]
$\D$          & unit disk in $\C$, i.e., $\D=K_1(0)$& \pageref{def:D} \\[1mm]
$\D'$          & punctured unit disk in $\C$, i.e., $\D'=\D \backslash \{ 0\}$& \pageref{def:D'} \\[1mm]
$\D_R$        & open disk in $\C$, center $0$, radius $R>0$& \pageref{def:DR}\\[1mm] 
$\D'_R$        & punctured disk $\{z \in \C \, : \, 0<|z|<R\}$ & \pageref{def:DR'}\\[1mm] 
$\Delta_R$    & complement of $\overline{\D}_R$, i.e., $\{z \in \C \, : |z|>R\}$& \pageref{def:DeltaR} \\[1mm]
$A_{r,R}$     & annulus $\{z \in \C \, : \, r<|z|<R\}$ & \pageref{def:ArR} \\[1mm]
$ \Delta$    & Laplace Operator $\Delta:=\frac{\partial^2}{\partial x^2}+\frac{\partial^2}{\partial y^2}$ & \pageref{def:laplace}\\
             & or  generalized Laplace Operator & \pageref{def:genlap}\\[1mm] $\kappa_{\lambda}$ & (Gauss) curvature of $\lambda(z) \, |dz|$ & \pageref{def:curvature} \\[1mm]
$f^*\lambda$ & pullback of $\lambda(z) \, |dz|$ under $f$ & \pageref{eq:pullback}\\[1mm]
$\lambda_G(z) \, |dz|$ & hyperbolic metric of $G$ & \pageref{def:hypi} \\[1mm]
$L_{\lambda}(\gamma)$ & $\lambda$--length of path $\gamma$ & \pageref{def:leng} \\[1mm]
$\d$ & distance function associated to $\lambda(z) \, |dz|$ & \pageref{def:dist}\\[1mm]
$C(M)$ & set of continuous functions on $M\subseteq \C$ & \pageref{def:C}\\[1mm]
$C^k(M)$ & set  of functions having all derivatives of order $\le k$ continuous & \pageref{ref:diff} \\ 
&     in $M \subseteq \C$ & \\[1mm]
$m{_\zeta}$ & two--dimensional Lebesgue measure w.r.t.~$\zeta \in \C$ & \pageref{def:lm}\\
\end{tabular}

\newpage

\section{Curvature} \label{sec:curvature}

In the sequel $D$ and $G$ \label{def:G} always denote  domains in the complex plane $\C$.

\subsection{Conformal metrics}

Conformal maps preserve angles between intersecting paths\footnote{Throughout, all paths 
are  assumed to be
continuous and piecewise continuously differentiable.}, but the
euclidean length
$$ \int \limits_{\gamma} \, |dz|:=\int \limits_{a}^{b} |\gamma'(t)| \, dt$$
of a path
$\gamma : [a,b] \to \C$ is in general not conformally invariant.
It is therefore advisable to allow more flexible ways to measure the length of paths.

\begin{definition}[Conformal densities]
A continuous function $\lambda : G \to [0,+\infty)$, $\lambda \not\equiv 0$, is
called conformal density on $G$.
\end{definition}

\begin{definition}[Metrics and Pseudo--metrics]
Let $\gamma : [a,b] \to G$ be a path in $G$ and $\lambda$
a conformal density on $G$. Then
$$ L_{\lambda}(\gamma):=\int \limits_{\gamma} \lambda(z) \, |dz|:=\int
\limits_{a}^{b} \lambda(\gamma(t)) \, |\gamma'(t)| \, dt$$
is called the $\lambda$--length of the path $\gamma$ and the quantity \label{def:leng}
$\lambda(z) \, |dz|$ is called  conformal pseudo--metric on $G$.
If $\lambda(z)>0$ throughout $G$, we say $\lambda(z) \, |dz|$ is a
conformal metric on $G$. We call a conformal
pseudo--metric $\lambda(z) \, |dz|$ on $G$ regular, if $\lambda$ is of class
$C^2$ in $\{z \in G \, : \, \lambda(z)>0\}$.
\end{definition}

\begin{remark1}
 The $\lambda$--length of $\gamma$ only depends on the trace of $\gamma$, but not on its parameterization.
\end{remark1}

Let $\lambda(w) \, |dw|$ be a conformal pseudo--metric
on $D$  and $w=f(z)$ a non--constant analytic map from $G$ to $D$.
We want to  find a conformal pseudo--metric $\mu(z) \, |dz|$ on $G$ such that
$L_{\mu}(\gamma)=L_{\lambda}(f \circ \gamma)$ for every path $\gamma$ in
$G$. There is clearly only one possible choice for $\mu(z) \, |dz|$ since by
the change--of--variables formula
$$ L_{\lambda}(f\circ \gamma)=\int \limits_{f \circ \gamma} \lambda(w) \,
|dw|=
\int \limits_{\gamma} \lambda(f(z)) \, |f'(z)| \, |dz|
=L_{\lambda\circ f \cdot |f'|}(\gamma) \, .$$

\begin{definition}[Pullback of conformal pseudo--metrics]
Let $\lambda(w) \, |dw|$ be a conformal pseudo--metric
on $D$  and $w=f(z)$ a non--constant analytic map from $G$ to $D$. Then
the conformal pseudo--metric
\begin{equation} \label{eq:pullback}
 (f^*\lambda)(z) \, |dz|:=\lambda(f(z)) \, |f'(z)| \, |dz|
\end{equation}
is called the pullback of $\lambda(w) \, |dw|$ under $w=f(z)$.
\end{definition}

Thus
\begin{equation} \label{eq:length}
 L_{f^*\lambda}(\gamma)=L_{\lambda}(f \circ\gamma)\, . 
\end{equation}

\subsection{Gauss curvature} 

Let $\lambda(w) \, |dw|$ be a conformal pseudo--metric.
 We wish to introduce  a quantity $T_{\lambda}$ which is conformally invariant  in the sense that
$$ T_{f^*\lambda}(z)=T_{\lambda}(f(z)) \, $$
for all conformal maps $w=f(z)$.

\medskip

Consider (\ref{eq:pullback}). We need to
 eliminate the conformal factor $|f'(z)|$. For this we note that $\log |f'|$
is harmonic, so $\Delta (\log |f'|)=0$ and therefore
$$ \Delta (\log f^*\lambda)(z)=\Delta (\log \lambda \circ f)(z)+\Delta (\log
|f'|)(z)=\Delta (\log \lambda \circ f)(z)=\Delta (\log \lambda) (f(z)) \, |f'(z)|^2 \, .
$$
In the last step we used the chain rule $\Delta (u\circ f)=(\Delta u \circ f) \cdot |f'|^2$
for the Laplace Operator $\Delta$ \label{def:laplace} and holomorphic functions $f$.
We see that $T_{\lambda}(z):=\Delta (\log\lambda)(z)/\lambda(z)^2$ is conformally invariant.

\begin{definition}[Gauss curvature]\label{def:curvature}
Let $\lambda (z) \, |dz|$ be a regular conformal pseudo--metric on $G$. Then
$$ \kappa_{\lambda}(z):=-\frac{\Delta (\log \lambda)(z)}{\lambda(z)^2}$$
is defined for all points $z \in G$ where $\lambda(z) > 0$. The quantity $\kappa_{\lambda}$  is
called the (Gauss) curvature of $\lambda(z) \, |dz|$.
\end{definition}

Our preliminary considerations can now be summarized as follows.

\begin{theorem}[Theorema Egregium]
For every analytic map $w=f(z)$ and every regular conformal pseudo--metric $\lambda(w)\,|dw|$ the relation
\begin{equation} \label{eq:curv}
\kappa_{f^*\lambda}(z)=\kappa_{\lambda}(f(z))
\end{equation}
is satisfied provided $\lambda(f(z)) \,  |f'(z)|>0$. 
\end{theorem}



For the {\it euclidean metric} $\lambda(z) \, |dz|:=|dz|$ one easily finds
$\kappa_{\lambda}(z) = 0$. 
Note that every regular conformal metric $\lambda(z) \, |dz|$
with {\it non--positive} curvature gives rise to a {\it subharmonic}
function $u(z):=\log \lambda(z)$ of class $C^2$ and vice versa.

\subsection{Constant curvature} \label{subsec:constcurv}

In view of (\ref{eq:curv}) conformal metrics $\lambda(z) \, |dz|$ with {\it constant}
curvature are of particular relevance, since for such metrics curvature is
an {\it absolute} conformal invariant, i.e., $\kappa_{f^*\lambda} \equiv
\kappa_{\lambda}$. In order to find all constantly curved conformal metrics, we
set  $u(z):=\log \lambda(z)$ and are therefore led to the problem 
of computing all solutions $u(z)$ to the equation
\begin{equation*} 
\Delta u=-k \, e^{2 u} \, ,
\end{equation*} 
where $k \in \R$ is a real \label{ref:R} constant. This equation  is called
{\it Liouville's equation}.

\begin{remark1}[Curvature zero]
If $k=0$, then the solutions of Liouville's equation 
 are precisely all harmonic functions. Thus, a regular conformal
metric $\lambda(z) \, |dz|$ has zero curvature if and only if $\lambda=e^u$
 for some harmonic function $u$.
\end{remark1}

We next consider the case of constant negative curvature $k<0$. We may
normalize and choose $k=-1$. Thus we need to consider
$$ \Delta u= e^{2 u} \, .$$
We shall later find all solutions to this equation. For now it suffices to 
 determine all {\it radially symmetric} solutions $u(z)=u(|z|)$, which can be
computed rather easily. As a byproduct we obtain in a constructive  way
 a number of important examples of conformal metrics.
These metrics will play an ubiquitous
r$\hat{\mbox{o}}$le in the sequel.

\begin{proposition}[Radially symmetric
  metrics with constant curvature {\mathversion{bold}$-1$}] \label{prop:1}
Let $\lambda(z) \, |dz|$ be a radially symmetric regular conformal metric
defined on some open annulus centered at $z=0$ with constant curvature $-1$. Then one of the following holds.
\begin{itemize}
\item[(a)]  $\lambda$ is defined on a disk $\D_R:=\{z \in \C \, : \, |z|<R\}$  and \label{def:DR}

\begin{tabular}{rl}  \hspace*{2cm} & $\displaystyle \lambda(z)=\lambda_{\D_R}(z):=\frac{2\, R}{R^2-|z|^2} \, .$ \\
\end{tabular}
\item[(b)] $\lambda$ is defined on a punctured disk $\D_{R}':=\{ z \in \C \,
  :   \, 0<|z|<R\}$  \label{def:DR'}and either

\begin{tabular}{rl}   \hspace*{2cm} & $\displaystyle \lambda(z) =\lambda_{\D'_R}(z):=\frac{1}{ |z| \log (R/|z|)}$ \\
\hspace*{-0.5cm} or  \hspace*{1.8cm} & \\  \hspace*{2cm} 
 &$\displaystyle  \lambda(z) =\frac{2\, \alpha R^{\alpha} |z|^{\alpha-1}}{R^{2 \alpha}- |z|^{2
    \alpha}} \,  \qquad  \text{ for some }\alpha \in (0, +\infty) \backslash \{ 1 \} \, .$
\end{tabular}
\item[(c)] $\lambda$ is defined on an annulus $A_{r,R}:=\{ z \in \C \, : \, r<|z|<R\}$ and \label{def:ArR}

\begin{tabular}{rl}   \hspace*{2cm}  & $\displaystyle \lambda(z)=\lambda_{A_{r,R}}(z):=\frac{\pi}{\log (R/r)} \frac{1}{|z| \, \sin \left[ \pi \displaystyle \frac{\log (R/|z|)}{\log (R/r)} \right]}   \, .$ \\
\end{tabular}
\item[(d)] $\lambda$ is defined on $\Delta_R:=\{ z \in \C \, : \, R<|z|<+\infty\}$ for some $R>0$ and either \label{def:DeltaR}

\begin{tabular}{rl}    \hspace*{2cm} & $ \displaystyle \lambda(z) =\lambda_{\Delta_R}(z):=\frac{1}{ |z| \log (|z|/R)}$ \\
\hspace*{-0.5cm} or  \hspace*{2cm} & \\  \hspace*{2cm} 
 &$\displaystyle  \lambda(z) =\frac{2\, \alpha R^{\alpha} |z|^{\alpha-1}}{|z|^{2 \alpha}- R^{2
    \alpha}} \,  \qquad  \text{ for some }\alpha \in (0, +\infty)  \, .$
\end{tabular}
\end{itemize}
\end{proposition}

\medskip

{\bf Proof.} We give the main steps involved in the proof and leave the
details to the reader. Since $\lambda$ is radially symmetric
we get for $u(z):=\log \lambda(z)$
$$ \Delta u= \frac{1}{r} \left( r u'(r) \right)' \, ,$$
where $r=|z|$ and $'$ indicates differentiation with respect to $r$. Thus
Liouville's equation $\Delta u=e^{2 u}$ transforms into the ODE
$$ (r u'(r))'= r \, e^{2 u(r)} \, .$$
We substitute $r=e^{x}$ and obtain for $v(x):=u(e^x)+x$ the ODE
$$ v''(x)=  \, e^{2 v(x)} \, .$$
This ODE has $2 \, v'(x)$ as an integrating factor,  so
$$ \left( v'(x)^2 \right)'= \left( e^{2 v(x)} \right)' \, .$$ 
Integration yields
$$ v'(x)=\pm \sqrt{ e^{2 \, v(x)}+c} $$
with some constant of integration $c \in \R$. The resulting two ODEs (one for
each sign) are separable and can be solved by elementary integration.
This leads to explicit formulas for $v(x)$ and thus also for
$\lambda(z)=e^{u(|z|)}=e^{v(\log |z|)}/|z|$. \hfill{$\blacksquare$}

\begin{definition}[The hyperbolic metric]
We call
\begin{equation*} \label{eq:hypdef}
\lambda_{\D} (z)\, |dz|:=\frac{2}{1-|z|^2}\, |dz|
\end{equation*}
the hyperbolic metric for the unit disk $\D:=\{z \in \C \, : \, |z|<1\}$ and 
$\lambda_{\D}$ \label{def:D} the hyperbolic density of $\D$. The curvature of
$\lambda_{\D}(z) \, |dz|$ is $-1$.
\end{definition}

\begin{warning}
 Some authors call 
$$ \frac{ |dz|}{1-|z|^2}$$
the hyperbolic metric of $\D$. This metric has constant curvature $-4$.
\end{warning}

\section*{Exercises for Section 2}
\begin{enumerate}
\item \label{exe:1}
 Let $T$ be a conformal self--map of $\D$. Show that $T^*\lambda_{\D}\equiv \lambda_{\D}$.
\item \label{exe:2} Denote \label{def:D'} by $\D':=\D \backslash \{ 0\}$ the punctured unit disk and let
$$ \lambda_{\D'}(z) \, |dz|:=\frac{|dz|}{ |z| \log(1/|z|)}
  \, . $$ 
 Find an analytic function $\pi : \D \to \D'$ such that
 $\pi^*\lambda_{\D'}=\lambda_{\D}$ and $\pi(0)=e^{-1}$.\\
(Consider the ``ODE'' $$\lambda_{\D'}(\pi(x))\, |\pi'(x)|= \lambda_{\D}(x)$$
  for $x \in (-1,1)$ assuming that $\pi(x)$ is real and positive  and $\pi'(x)$ is real and negative.)

\item  \label{exe:2.3} Find all radially symmetric regular conformal metrics 
 defined on some open annulus centered at $z=0$ with constant curvature $+1$.

\item \label{exe:pickeq}
Let $f : \D \to \C$ be an analytic map with $f(0)=0$ and
$f^*\lambda_{\D} \equiv \lambda_{\D}$ in a neighborhood of $z=0$.
Show that $f(z)=\eta z$ for some $|\eta|=1$.
\item \label{exe:geodesic_curvature} (Geodesic curvature)

Let $\lambda(z) \, |dz|$ be a regular conformal metric on $D$ and $\gamma : [a,b] \to D$ a $C^2$--path. Then
$$ \kappa_{\lambda}(t,\gamma):=\frac{\Im \left[ \frac{d}{dt} \big( \log \gamma'(t) \big)+2 \frac{\partial \log \lambda}{\partial z}(\gamma(t)) \cdot \gamma'(t) \right]}{\lambda(\gamma(t)) \, |\gamma'(t)|} \, , \qquad t \in [a,b] \, , $$
is called the geodesic curvature of $\gamma$ at $t$ for $\lambda(z) \, |dz|$.
Show that if $f$ is an analytic map, then
$\kappa_{f^*\lambda}(t,\gamma)=\kappa_{\lambda}(t,f \circ \gamma)$.
\end{enumerate}

\section{The Fundamental Theorem} \label{sec:ahlfors}

\subsection{Ahlfors' Lemma}

The hyperbolic metric $\lambda_{\D}(z) \, |dz|$ has an important extremal
property.

\begin{lemma}[Ahlfors' Lemma \cite{Ahl38}] \label{lem:ahlfors}
Let $\lambda(z) \, |dz|$ be a regular  conformal pseudo--metric on $\D$ with
curvature bounded above by $-1$. Then $\lambda(z) \le\lambda_{\D}(z)$ for every $z \in \mathbb{D}$.
\end{lemma}

Following A.~Beardon and D.~Minda we call Lemma \ref{lem:ahlfors} the {\it
  Fundamental Theorem}.

\medskip


{\bf Proof.}
Fix $0<R<1$ and consider $\lambda_{\D_R}(z)= 2R/(R^2-|z|^2)$ on the disk $\D_R$.
 Then the  function
$u : \D_R \to \R \cup \{-\infty\}$ defined by
$$u(z):=\log \left( \frac{\lambda(z)}{\lambda_{\D_R}(z)} \right) \, , \qquad z \in \D_R \, , $$

tends to $-\infty$ as  $|z|\to R$ and therefore attains its maximal value
at some point $z_0 \in \D_R$ where $\lambda(z_0)>0$.
 Since 
$u$ is of class $C^2$ in a neighborhood of $z_0$, we get
$$ 0 \ge \Delta u(z_0)=\Delta \log \lambda(z_0)-\Delta \log
\lambda_{\D_R}(z_0) \ge   \lambda(z_0)^2-\, \lambda_{\D_R}(z_0)^2 \, .$$
Thus  $u(z) \le u(z_0)\le 0$ for all $z\in \D_R$, and so
$\lambda(z) \le \lambda_{\D_R}(z)$ for $z \in \D_R$. 
Now, letting $R \nearrow 1$, gives $\lambda(z) \le \lambda_{\D}(z)$ for all $z \in \D$. \hfill{$\blacksquare$}


\begin{lemma} \label{lem:equality}
Let $\lambda(z) \, |dz|$ be a regular conformal pseudo--metric on $\D$ with
curvature $\le -1$ and $\mu(z) \, |dz|$ a regular conformal metric on $\D$ with
curvature $=-1$ such that $\lambda(z) \le \mu(z)$ for all $z \in \D$. Then
 either $\lambda<\mu$ or $\lambda \equiv \mu$.
\end{lemma}

In particular, if equality holds in Lemma \ref{lem:ahlfors} for one point $z \in \D$, then $\lambda \equiv \lambda_{\D}$.

\medskip

{\bf Proof.}
Assume $\lambda(z_0)=\mu(z_0)>0$ for some point $z_0 \in \D$. We may take
$z_0=0$, since otherwise we could consider $\tilde{\lambda}:=T^*\lambda$ and
$\tilde{\mu}:=T^*\mu$ for a conformal self--map $T$ of $\D$ with $T(0)=z_0$.
Let $u(z):=\log (\mu(z)/\lambda(z))\ge 0$ and 
$$ v(r):=\frac{1}{2 \pi} \int \limits_{0}^{2 \pi} u(r e^{i t}) \, dt \,
, \qquad V(r):=\int \limits_{0}^r v(\rho) \, d\rho \, .$$
Then $v(0)=0 \le v(r)=V'(r)$. To show $V \equiv 0$ observe
$ \Delta u\le \mu(z)^2 \left(1- e^{-2 u} \right) \le C u$
in $|z| <R<1$ where $C:=2 \sup \{\mu(z)^2 \, : \, |z| < R\}<+\infty$
for small $R>0$. Hence an application of
Green's theorem, \label{def:lm}
$$r v'(r)=\frac{r}{2 \pi} \frac{d}{dr} \int \limits_{0}^{2 \pi} u(r e^{it}) \, dt=
\frac{1}{2 \pi} \iint \limits_{|\zeta|<r}
\Delta u(\zeta) \, dm_{\zeta}=
\frac{1}{2 \pi} \int \limits_{0}^r \int \limits_{0}^{2 \pi} \Delta u(\rho e^{it}) \, dt
\, \rho \, d\rho\, ,  $$
 leads for all $0 \le r < R$ to the estimate
$$
r V''(r)-C r V(r) \le 
r v'(r)-C \int \limits_{0}^r \rho \, v(\rho) \, d\rho=\frac{1}{2\pi} 
\int \limits_{0}^{r} \int \limits_{0}^{2 \pi} \left( \Delta u(\rho e^{i t})-C u(\rho
e^{i t}) \right) \, dt \, \rho
\, d\rho \le 0 \, .$$
Thus $V''(r) \le C\,  V(r)$, so $(V'(r)^2)' \le C (V(r)^2)'$. This implies
$V'(r) \le \sqrt{C} \cdot V(r)$, i.e., 
$(V(r) e^{-\sqrt{C}\,  r})' \le 0$ and therefore $0 \le
V(r)e^{-\sqrt{C} \, r} \le V(0) =0$. Hence $V \equiv 0$ and $v \equiv 0$ on $[0,R)$.
Thus  $u \equiv 0$ first on $|z|< R$ and therefore in all of $\D$, so $\lambda
\equiv \mu$. 
~\hfill{$\blacksquare$}

\subsection{Applications of Ahlfors' Lemma}

\begin{corollary}[Pick's Theorem]\label{cor:pick}
Let $f : \D \to \D$ be an analytic function. Then
$$ \frac{|f'(z)|}{1-|f(z)|^2} \le \frac{1}{1-|z|^2} \, , \quad z \in \mathbb{D}\, .$$
\end{corollary}

{\bf Proof.}
 Note that for a non--constant analytic function $f$,  $(f^*\lambda_{\D})(z) \, |dz|$ is
 a regular conformal pseudo--metric on $\D$
with constant curvature $-1$. Now apply the Fundamental Theorem.~\hfill{$\blacksquare$}



\begin{corollary} \label{thm:nometriconC}
The complex plane $\C$ admits no regular conformal pseudo--metric 
with curvature  $\le -1$.
\end{corollary}

{\bf Proof.} Assume that there is  a regular conformal pseudo--metric 
$\lambda(z) \, |dz|$ 
with curvature  $\le -1$ on $\C$. For fixed $R>0$ consider $f_R(z):=R z$ for $z \in \D$. Then $(f_R^*\lambda)(z)
\, |dz|$ is a regular conformal pseudo--metric on $\D$ with curvature $\le
-1$. The Fundamental Theorem gives
 $\lambda(Rz) \,  R =(f^*_R\lambda)(z) 
\le \lambda_{\D}(z)$ for every $z \in
\D$, that is,  $\lambda(w)  \le \lambda_{\D}(w/R)/R$ for each $|w|<R$. 
Letting  $R \to +\infty$  we thus obtain $\lambda \equiv 0$, which is not possible.
\hfill{$\blacksquare$}

\begin{corollary}\label{thm:liouvilleabstract}
Assume $G$ admits a regular conformal metric with
curvature bounded above by $-1$. Then
every entire function $f : \C\to G$ is constant.
\end{corollary}

{\bf Proof.} Let $\lambda(w) \, |dw|$ be a regular conformal metric on $G$ with
curvature bounded above by $-1$. If there would be 
a non--constant  analytic function  $f : \C \to G$, then 
$(f^*\lambda)(z) \, |dz|$ would be a regular conformal pseudo--metric on $\C$
with curvature  $\le -1$, contradicting Corollary \ref{thm:nometriconC}.
\hfill{$\blacksquare$}

\medskip

Thus the exponential map shows that there is no regular conformal metric on 
the punctured plane $\C':=\C \backslash \{ 0 \}$ \label{def:C'} with
curvature bounded above by $-1$. However, there is such a metric on the twice--punctured plane.

\begin{theorem} \label{thm:hyp}
The twice--punctured plane $\C'':=\mathbb{C}\backslash\{0,1\}$ carries a
regular conformal metric with curvature $\le -1$. \label{def:C''}
In particular, every domain $G \subset \C$ with at least two boundary points carries a
regular conformal metric with curvature $\le -1$. 
\end{theorem}

{\bf Proof.}
We first show  that for  $\eps>0$  sufficiently small 
$$ \tau(z):= \eps \, \frac{\sqrt{1+ \left| z \right| ^{\frac{1}{3}}}}{ \left| z \right| ^{\frac{5}{6}}} \frac{
\sqrt{1+ \left| z-1 \right| ^{\frac{1}{3}}}}{\left| z-1 \right|^{\frac{5}{6}}}
$$
defines a regular conformal metric $\tau(z) \, |dz|$ 
 with curvature $\le -1$ on $\C''$.

\smallskip

A quick computation using polar coordinates gives for
 $r:= \left| z \right|$ 
\begin{eqnarray*} \Delta \, \log \frac{\sqrt{1+ \left| z \right|^{\frac{1}{3}}}}{\left| z \right| ^{\frac{5}{6}}} &=& \Delta \left( \frac{1}{2}\log
 \left( 1+r^{\frac{1}{3}} \right) - \frac{5}{6} \log r \right)= 
 \frac{1}{2} \Delta \left[ \log \left( 1+r^{\frac{1}{3}} \right) \right] = \\
&=& \frac{1}{2} \left( \frac{\partial^2}{ \partial r^2}
 \log (1+r^{\frac{1}{3}})  + \frac{1}{r}
\frac{\partial}{\partial r} \log(1+r^{\frac{1}{3}}) \right)= 
\frac{1}{18} \frac{1}{(1+|z|^{\frac{1}{3}})^2 |z|^{\frac{5}{3}}}\, 
\end{eqnarray*}
and analogously
\begin{eqnarray*} \Delta \, \log \frac{\sqrt{1+ \left| z-1 \right|^{\frac{1}{3}}}}{\left| z-1 \right| ^{\frac{5}{6}}} &=&
\frac{1}{18} \frac{1}{(1+|z-1|^{\frac{1}{3}})^2 |z-1|^{\frac{5}{3}}}\, .
\end{eqnarray*}
Thus
$$ \kappa_{\tau}(z)= - \frac{1}{18 \, \eps^2}  \left[ 
\frac{|z-1|^{5/3}}{(1+|z|^{\frac{1}{3}})^3 
(1+|z-1|^{\frac{1}{3}})}+ \frac{|z|^{5/3}}{(1+|z|^{\frac{1}{3}})
(1+ |z-1|^{\frac{1}{3}})^3} \right]
\le -1$$
for each $z \in \C''$, if $\eps>0$ is small enough since
$$ \lim \limits_{z \to 0}  \kappa_{\tau}(z)=-\frac{1}{36  \cdot\eps^2}, \qquad
\lim \limits_{z \to 1}  \kappa_{\tau}(z)=-\frac{1}{36 \cdot\eps^2}, \qquad
\lim \limits_{|z| \to +\infty}  \kappa_{\tau}(z)=-\infty.$$
If $a,b \in \C$ with $a\not=b$, then the pullback of $\tau(w) \, |dw|$ under
$w=T(z)=(z-a)/(b-a)$ yields a regular conformal metric with curvature $\le -1$
on $\C\backslash \{ a,b\}$.
Thus every domain $G$ with at least two boundary points carries 
a regular conformal metric with curvature $\le -1$ on $G$.~\hfill{$\blacksquare$}

\begin{corollary}[The Little Picard  Theorem]
Every entire function $f : \C \to \C$ which omits two distinct complex numbers is constant.
\end{corollary}
{\bf Proof.} Let $a$ and $b$ be two distinct complex numbers and
 assume that $f(z)\not=a$ and $f(z) \not= b$ for every $z \in \C$. 
By Theorem \ref{thm:hyp}
$\C \setminus \{ a , b \}$ carries a regular conformal metric with curvature
$\le -1$, so  $f$ must be
constant in view of Corollary
\ref{thm:liouvilleabstract}. \hfill{$\blacksquare$}

\begin{corollary}[Huber's Theorem] \label{huber}
Let $G \subseteq \C$ be a domain which carries a  regular conformal metric with
curvature $\le -1$ and let 
$f : \D' \to
G$ be a holomorphic function. If there is a sequence $(z_n ) \subset
\D'$ such that $\lim_{n \to \infty} z_n=0$ and such that 
$\lim_{n \to \infty} f(z_n) $ exists and belongs to $G$, 
then $z=0$ is a removable singularity of $f$.
\end{corollary}
{\bf Proof.} By hypothesis $G$ carries a 
   conformal metric  $\lambda(w) \, |dw|$ on $G$ with curvature $\le -1$.
We may assume w.l.o.g.~that $r_n:=|z_n|$ is  monotonically decreasing to $0$ 
 and that $f(z_n) \to 0 \in G$. We consider the closed curves
 $\gamma_n:=f(\partial \D_{r_n})$ in $G$. \label{def:p}
Since $(f^*\lambda)(z) \, |dz|$ is a regular conformal pseudo--metric on $\D'$ with curvature $\le -1$, we get (see Exercise \ref{sec:ahlfors}.\ref{exe:punc}) 
$$(f^*\lambda) (z) \le \lambda_{\D'}(z)=\frac{1}{|z| \log(1/|z|)} \, , \qquad z \in \D' \, , $$
i.e., 
$$ L_{\lambda}(\gamma_n) =L_{f^*\lambda}(\partial \D_{r_n})
\le L_{\lambda_{\D'}}(\partial \D_{r_n}) =
\frac{2 \pi}{\log(1/r_n)} \to 0 \qquad \text{ as } n \to  \infty \, . $$
Now let $K$ be a closed disk of radius $\eps>0$ centered at $w=0$ which is compactly contained in
$G$. We may assume $f(z_n) \in K$ for any $n \in \N$. \label{def:N} Since $\lambda(w) \, |dw|$ is a
conformal {\it metric}, its density is bounded away from zero on $K$, so $\lambda(w) \ge
c>0$ for all $w \in K$. Hence the euclidean length of $\gamma_n \cap K$
is bounded above by $L_{\lambda}(\gamma_n)/c$. 
In particular, $\gamma_n=f(\partial \D_{r_n}) \subset K$, i.e., $|f(z)| \le \eps$ on $|z|=r_n$
 for all but finitely many indices $n$. The maximum principle implies that
$|f(z)| \le \eps$ in a punctured neighborhood of $z=0$, so the singularity
$z=0$ is removable.~\hfill{$\blacksquare$}

\begin{corollary}[The Big Picard Theorem]
If an analytic function $f : \D' \to \C$ has an essential singularity at $z=0$, 
then there exists at most one complex number $a$ such that the equation
$f(z)=a$ has no solution. 
\end{corollary}

{\bf Proof.}
Assume $a,b \not \in f(\D')$ for $a\not=b$. If $c \in \C \backslash \{a,b\}$,  then the theorem of Casorati--Weierstra{\ss} tells us that there
exists a sequence $( z_n)$ such that $z_n \to 0$  and $f(z_n) \to c$. 
This, however,  contradicts Huber's Theorem (Corollary \ref{huber}), because $\C \setminus \{ a, b \}$
carries a  regular conformal metric with curvature $\le -1$  by Theorem  \ref{thm:hyp}.
\hfill{$\blacksquare$}

\section*{Exercises for Section 3}

\begin{enumerate}
\item \label{exe:punc} 

The Fundamental Theorem clearly holds for every disk $\D_R$, that is,
$\lambda(z) \le \lambda_{\D_R}(z)$ for every regular conformal pseudo--metric
$\lambda(z) \, |dz|$ with curvature $\le -1$ on $\D_R$.
\smallskip

\, Show that $\lambda(z) \le \lambda_{A_{r,R}}(z)$ for every regular conformal pseudo--metric
$\lambda(z) \, |dz|$ with curvature $\le -1$
 on the annulus $A_{r,R}$, see Proposition \ref{prop:1}. Deduce   that $\lambda(z) \le \lambda_{\D'}(z)$ for every regular conformal pseudo--metric
  $\lambda(z) \, |dz|$ with curvature $\le -1$ on $\D'$.

\item \label{exe:3.2} Let $\lambda(z) \, |dz|$ be a regular conformal metric with curvature $\le -1$ on a disk $K_R(z_0):=\{z \in \C \, : \, |z-z_0|<R\}$ \label{def:dik} with radius $R>$0 and center $z_0 \in \C$. Define $u(z):=\log \lambda(z)$ and
$$ v(r):=\frac{1}{2 \pi} \int \limits_{0}^{2 \pi} u\left(z_0+r e^{i t} \right) \, dt \, .$$
Show that 
$$ v(r)-v(0) \ge \frac{e^{2 v(0)}}{4}\,   r^2 \,  \qquad \text{ for } 0 \le r < R\,.$$
(Hint: Compute $(r v'(r))'/r$ and recall Jensen's inequality.)

\item  \label{exe:3.3} Modify the proof of  Lemma \ref{lem:equality} to prove the 
following special case of the maximum
  principle of E.~Hopf \cite{Hopf} (see also Minda \cite{Min87}):

{\sl Let $c : G \to \R$ be a non--negative bounded function and $u : G \to \R$
  of class $C^2$ such that
$\Delta u \ge c(z) u$
  in $G$. If $u$ attains a non--negative maximum at some interior point of $G$,
  then $u$ is constant.}

\item \label{exe:3.4} Use Lemma \ref{lem:equality} and Exercise
  \ref{sec:curvature}.\ref{exe:pickeq} to describe all cases of equality
in Pick's Theorem (Corollary \ref{cor:pick}) for some $z \in \D$.

\item \label{exe:robinson}
Let $z_1, \ldots , z_{n-1}$ be distinct points in $\C$ and $z_n=\infty$. Further,
let $\alpha_1, \ldots \, ,\alpha_{n}$ be real numbers $< 1$ with $s:=\alpha_1+\cdots +\alpha_n>2$. 
\begin{itemize}
\item[(a)]
Define for $\eps>0$ and $\delta>0$
$$ \tau(z):=\eps \prod \limits_{j=1}^{n-1}
\frac{\left[ 1+|z-z_j|^{\delta} \right]^{\frac{s-2}{(n-1) \delta}}}{|z-z_j|^{\alpha_j}} \, .$$
Show that if $\eps$ and $\delta$ are sufficiently small, then $\tau(z) \, |dz|$
is a regular conformal metric with curvature $\le -1$
on $\C \backslash \{z_1, \ldots, z_{n-1}\}$ such that
$|z-z_j|^{\alpha_j} \tau(z)$ is continuous and positive at $z=z_j$ for each $j=1, \ldots , n-1$ and $|z|^{2-\alpha_n} \tau(z)$ is continuous and positive at $z=\infty$.
\item[(b)] Let $f$ be a meromorphic function on $\C$. We say that $w$ is an exceptional value of order $m$ if the equation $f(z)=w$ has no root of multiplicity less than $m$. Prove the following result of Nevanlinna:
{\it Let $f$ be a meromorphic function on $\C$ with exceptional values
$w_j$ of order $m_j$. If 
$$ \sum \limits_{j} \left( 1-\frac{1}{m_j} \right)>2 \, , $$
then $f$ is constant.}
\end{itemize}
\end{enumerate}

\section*{Notes}
Lemma \ref{lem:ahlfors} was established by Ahlfors \cite{Ahl38} in 1938.
He used it to derive quantitative bounds in the theorems of Bloch and Schottky. 
Lemma \ref{lem:equality} is due to M.~Heins \cite{Hei62}. The proof  given here is a simplified version of Heins' method. Different proofs were also found by D.~Minda \cite{Min87} and
H.~Royden \cite{Roy86}. An elegant direct treatment of the case of equality in Ahlfors' lemma was given by H.~Chen \cite{Ch01}. However, his method
does not appear to be strong enough to prove the more general Lemma \ref{lem:equality}. The proof of Theorem \ref{thm:hyp} is from the paper \cite{MinScho3}
of D.~Minda and G.~Schober, which in turn is based on the work of R.~Robinson \cite{Rob39}. The statement of  Exercise \ref{sec:ahlfors}.\ref{exe:robinson} which generalizes Theorem \ref{thm:hyp}  is also due to Robinson.
Huber's Theorem (Corollary \ref{huber}) was originally proved by using the Uniformization Theorem, see \cite{Hub}. The proof here is essentially that of 
M.~Kwack \cite{Kwa}. 

\section{The hyperbolic metric} \label{sec:hyperbolic}

\subsection{SK--Metrics}

In many applications the class of  {\it regular} conformal metrics is too
restrictive to be useful. In analogy to the situation for subharmonic functions,
we consider for a real--valued function $u$  the integral means
$$ v(r):=\frac{1}{2\pi} \int \limits_{0}^{2 \pi} u(z_0+r e^{it}) \, dt \, .$$
Assuming momentarily that $u$ is of class $C^2$ an application of Green's
theorem gives 
$$ r v'(r)= \frac{1}{2 \pi} \iint \limits_{|\zeta|<r}
\Delta u(z_0+\zeta) \, dm_{\zeta} =
\frac{1}{2 \pi} \int \limits_{0}^r \int \limits_{0}^{2 \pi} \Delta u(z_0+\rho e^{it}) \, dt
\, \rho \, d\rho  \, . $$
We differentiate this identity w.r.t.~$r$ and get
$$ \Delta u(z_0)=\lim \limits_{r \to 0} \frac{1}{2 \pi} \int \limits_{0}^{2\pi} \Delta u(z_0+r
e^{it}) \, dt=
\lim \limits_{r \to 0} \frac{1}{r} \left( r v'(r)
  \right)'\,, $$
so $(r v'(r))'=r \Delta u(z_0)+o(r)$. We integrate from $0$ to $r$ and obtain
$v'(r)=(r/2) \Delta u(z_0)+o(r)$. Another integration from $0$ to $r$ yields
$v(r)=u(z_0)+(r^2/4) \Delta u(z_0)+o(r^2)$ and therefore
$$ \Delta u(z_0)=\lim \limits_{r \to 0} \frac{4}{r^2} \left( \frac{1}{2 \pi} \int_0^{2 \pi } u(z_0 +r e^{it}) \, dt -u(z_0)
\right) \, . $$
If $u$ is not $C^2$, then this limit need not exist. We thus define:

\begin{definition}\label{def:genlap}
Let $u : G \to \R$ be a
continuous function. Then the generalized lower Laplace Operator
$\Delta u$ of $u$ at a point $z \in G$ is defined by
$$\Delta u(z)= \liminf_{ r \to 0} \frac{4}{r^2} \left(
  \frac{1}{2 \pi} \int_0^{2 \pi } u(z +r e^{it}) \, dt - u(z)    
 \right) \, .$$
\end{definition}

\begin{remark1} 
Our preliminary considerations show that if $u$ is $C^2$ in a
neighborhood of a point $z \in G$,
then the generalized lower Laplace Operator coincides with the standard Laplace Operator at $z$. 
This justifies the use of the same symbol $\Delta$ for both operators.
\end{remark1}


\begin{definition} \label{def:gencurvature}
Let  $\lambda(z) \, |dz|$  be a  conformal pseudo--metric on
$G$. Then the (Gauss) curvature of $\lambda(z) \, |dz|$ at a point
$z \in G$ where $\lambda(z)>0$ is defined by
$$\kappa_{\lambda}(z)= - \frac{\Delta (\log
  \lambda)(z)}{\lambda(z)^2}\, .$$  
\end{definition}

For regular conformal metrics, Definition \ref{def:gencurvature}
is consistent with Definition \ref{def:curvature}.

\begin{definition}
A conformal pseudo--metric $\lambda(z) \, |dz|$ on $G$ is called SK--metric\footnote{``SK is intended to convey curvature subordinate to $-1$'', see \cite{Hei62}.} on
$G$, if its curvature is bounded above by $-1$.
\end{definition}

We now have all the technology to generalize the Fundamental Theorem
to SK--metrics.

\begin{theorem}[Fundamental Theorem] \label{lem:fund}
Let $\lambda(z) \, |dz|$ be an SK--metric on $\D$. Then $\lambda(z) \le\lambda_{\D}(z)$ for every $z\in \mathbb{D}$. In particular, $\lambda_{\D}(z)\, |dz|$ is the (unique) maximal SK--metric on $\D$.
\end{theorem}

Theorem \ref{lem:fund} can be proved in exactly the same way as Lemma
\ref{lem:ahlfors} by noting that for a continuous function $u$ the generalized
lower Laplace Operator is always non--positive at a local maximum. One could also
consider for $0<R<1$ the auxiliary function
$$ v(z):=\max \left\{ 0,\log \left( \frac{\lambda(z)}{\lambda_{\D_R}(z)} \right) \right\}\, .$$
The hypotheses guarantee that $v$ is {\it subharmonic} on $\D_R$ and $v=0$ on
$|z|=R$. By the maximum principle for subharmonic functions we get $v \le 0$
in $\D_R$, so $\lambda(z) \le \lambda_{\D_R}(z)$ there. Now, let $R \to 1$ in order to arrive at the desired result.

\subsection{The Perron method for SK--metrics}

We first discuss a number of simple, but very useful
techniques for producing new SK--metrics from old ones.

\begin{lemma} \label{lem:glue1}
Let $\lambda(z) \, |dz|$ and $\mu(z) \, |dz|$ be SK--metrics on $G$. Then
$\sigma(z):=  \max \{\lambda(z), \mu(z)\}$
induces   an SK--metric $\sigma(z) \, |dz|$ on $G$.
\end{lemma}

{\bf Proof.}
Clearly, $\sigma$ is continuous in $G$.
 If $\sigma(z_0)=\lambda(z_0)$ for a point $z_0 \in G$, then
\begin{eqnarray*}
\Delta\log \sigma(z_0)&=&  \liminf_{ r \to 0} \frac{4}{ r^2} \left(
  \frac{1}{2 \pi} \int_0^{2 \pi } \log \sigma(z_0 +r e^{it}) \, d\!\; t - \log
  \sigma(z_0)    
 \right)\\[2mm]
 & \ge &  \liminf_{ r \to 0} \frac{4}{ r^2} \left(
  \frac{1}{2 \pi} \int_0^{2 \pi } \log \lambda(z_0 +r e^{it}) \, d\!\; t - \log
  \lambda(z_0)    
 \right) \\[2mm] & \ge &   \, \lambda(z_0)^2 = \, \sigma(z_0)^2 \, .
\end{eqnarray*}
Thus $\kappa_{\sigma}(z_0) \le -1$. Similarly,
$\kappa_{\sigma}(z_0) \le -1$, if $\sigma(z_0)=\mu(z_0)$, so
in either case  $\kappa_{\sigma}(z_0) \le -1$.~\hfill{$\blacksquare$}

\medskip

From Lemma \ref{lem:glue1} it follows that SK--metrics need not be smooth.

\begin{lemma}[Gluing Lemma] \label{lem:glueing}
Let $\lambda(z)\, |dz|$  be an SK--metric  on  $G$
 and let $\mu(z)\, |dz|$ be an SK--metric on  an open subset $U$ of $G$ such that
the ``gluing condition''
\begin{equation*} 
\limsup \limits_{ U \ni z \to \xi} \mu(z) \le  \lambda(\xi)  
\end{equation*}
holds for all $\xi \in \partial U \cap G$. Then $\sigma(z)\, |dz|$ defined by
$$\sigma(z):=\begin{cases} \,   \max \{\lambda(z), \mu(z)\}   & \hspace{3mm} \, \text{for }  z \in  U \, , \\[2mm]
                       \,        \lambda(z)         & \hspace{3mm} \, \text{for } z \in G \backslash U
          \end{cases} $$
is an SK--metric on $G$.
\end{lemma}

{\bf Proof.}
We first note that the gluing condition guarantees that
 $\sigma$ is  continuous on $G$, so $\sigma(z) \, |dz|$ is a
 conformal pseudo--metric on $G$.

\medskip

We need to compute the curvature $\kappa_{\sigma}(z_0)$ of $\sigma(z) \, |dz|$ 
at each point $z_0 \in G$ with $\sigma(z_0)>0$.
If  $z_0 \in U$, then $\kappa_{\sigma}(z_0) \le -1 $ by Lemma \ref{lem:glue1}.
If   $z_0$ is an interior point of $G \backslash U$, then $\sigma(z)=\lambda(z)$ in
a neighborhood of $z_0$, so $\kappa_{\sigma}(z_0)=\kappa_{\lambda}(z_0) \le -1$.
Finally, if $z_0 \in \partial U \cap G$, then $\sigma(z_0)=\lambda(z_0)$ and thus
\begin{eqnarray*}
\Delta\log \sigma(z_0)&= &  \liminf_{ r \to 0} \frac{4}{ r^2} \left(
  \frac{1}{2 \pi} \int_0^{2 \pi } \log \sigma(z_0 +r e^{it}) \, d\!\; t - \log
  \sigma(z_0)    
 \right)\\[2mm]
 & \ge & \liminf_{ r \to 0} \frac{4}{ r^2} \left(
  \frac{1}{2 \pi} \int_0^{2 \pi } \log \lambda(z_0 +r e^{it}) \, d\!\; t - \log
  \lambda(z_0)    
 \right) \ge   \lambda(z_0)^2 = \sigma(z_0)^2
\end{eqnarray*}
Consequently,  $\kappa_{\sigma}(z_0) \le -1$.
\hfill{$\blacksquare$}

\medskip

By the Fundamental Theorem 
the hyperbolic metric $\lambda_{\D}(z) \, |dz|$ on the unit disk $\D$ is the
maximal SK--metric on $\D$. We shall next show that
 every domain $G$ which carries at least one SK--metric  carries 
 a maximal SK--metric $\lambda_G(z) \, |dz|$ as well. Furthermore, this maximal SK--metric
is a regular conformal metric with constant curvature $-1$.

\medskip

The construction hinges on a modification of Perron's method for subharmonic functions.
An important ingredient of Perron's method for subharmonic functions is played by
the Poisson integral formula, which allows  the construction of harmonic functions
with specified boundary  values on disks. We begin with the corresponding statement for
regular conformal metrics with constant curvature $-1$.

\begin{theorem} \label{thm:existence}
Let $K \subset  \C$ be an open  disk and $\varphi :\partial K \to (0, + \infty)$ a continuous function. Then there exists a unique regular conformal metric $\lambda(z)\, |dz|$ with constant curvature $-1$ on $K$ such that $\lambda$ is continuous on the closure $\overline{K}$ of $K$ \label{def:closure} and $\lambda(\xi)=\varphi(\xi)$ for all $\xi \in \partial K$. 
\end{theorem}

For the somewhat involved proof of Theorem \ref{thm:existence} we refer the reader to the Appendix.

\medskip

We next consider increasing sequences
of conformal metrics with constant curvature $-1$ and prove the following
theorem of Harnack--type.

\begin{lemma}\label{lem:monotone}
Let $\lambda_j(z)\, |dz|$, $j \in \mathbb{N}$, be a monotonically increasing
sequence of regular conformal metrics with constant curvature $-1$ on
$G$. Then $\lambda(z)\, |dz|$ for $\lambda(z):=\lim_{j \to \infty}
\lambda_j(z)$ is  a regular conformal metric of constant curvature $-1$ on $G$.  \end{lemma}

{\bf Proof.} 
 Fix $z_0 \in G$ and choose an open disk $K:=K_R(z_0)$ which is compactly
 contained in
$G$. By considering $T^*\lambda_j$ for  $T(z)=R \,z+z_0$ we may assume 
$K=\D$ and $\D_{R'} \subseteq G$ for some $R'>1$.
 The Fundamental Theorem (see Exercise \ref{sec:ahlfors}.\ref{exe:punc})
shows  $\lambda_j(z) \le \lambda_{\D_{R'}}(z)\le M<+\infty$ in $\D$, so
  $\lambda(z):=\lim_{j \to
  \infty} \lambda_j(z)$ is well--defined and bounded in $\D$.
We now consider to each $\lambda_j$ the function $u_j:= \log
\lambda_j$. Then by
Remark \ref{rem:representation}, 
\begin{equation*} 
u_j(z)=h_j(z)- \frac{1}{2\pi} \iint \limits_{\D} g(z, \zeta) \, e^{2 u_j(\zeta)}\,
dm_{\zeta}\, , \quad z \in \D \, ,
\end{equation*}
where $h_j:\D \to \R$ is harmonic in $\D$ and  continuous on $\overline{\D}$ with
$h_j(\xi)= u_j(\xi)$ for $\xi \in \partial \D$ and
 $g$ denotes Green's function of $\D$.
The fact that $(u_j)$ is monotonically increasing and bounded above
 implies that  $(h_j)$ forms a monotonically increasing
 sequence of harmonic functions bounded from above. Thus $(h_j)$ converges locally
uniformly in $\D$ to a harmonic function $h$. Letting $j \to \infty$  we
obtain by Lebesgue's Theorem on monotone convergence
$$u(z)=h(z)- \frac{1}{2\pi} \iint \limits_{\D} g(z, \zeta) \, e^{2 u(\zeta)}\,
dm_{\zeta}\, , \quad z \in \D\, $$
for $u(z)=\log \lambda(z)$.
It follows from Remark \ref{rem:representation} that $u$ is a $C^2$ solution to $\Delta u=e^{2u}$
in $\D$. Thus $\lambda(z) \, |dz|$ is  a regular conformal metric of constant
curvature $-1$ on $G$.~\hfill{$\blacksquare$}

\medskip

We finally prove an analog of the Poisson Modification of subharmonic
functions for SK--metrics.

\begin{lemma}[Modification]\label{lem:mod}
Let $\lambda(z)\, |dz|$ be an SK--metric on $G$ and $K$  an open disk which is compactly contained in $G$. Then there exists a unique SK--metric $M_K\lambda(z)\, |dz|$ on $G$ with the following properties:
\begin{itemize}
\item[(i)] $M_K\lambda(z)  =\lambda(z)$  for every $z \in G\backslash K$
\item[(ii)] $M_K\lambda(z)\, |dz| $ is a regular conformal metric on $K$ with constant curvature $-1$. 
\end{itemize}
\end{lemma}

We call $M_K\lambda$ the modification of $\lambda$ on $K$.\\[3mm]
\begin{figure}[h]
\centering
\subfigure[An SK--metric $\ldots$]{\includegraphics[width=7cm,height=7cm]{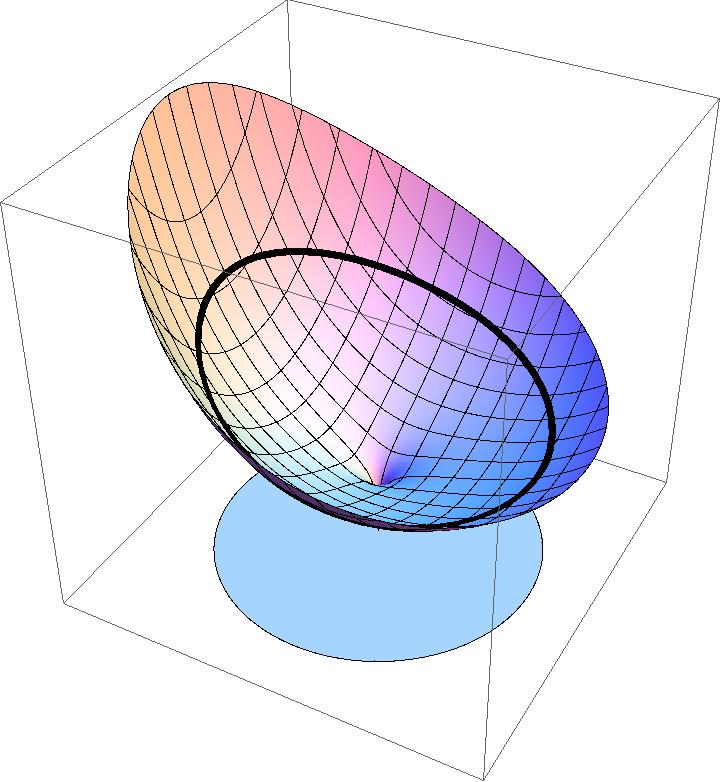}} \hspace{1cm}
\subfigure[$\ldots$ and its
modification]{\includegraphics[width=7cm,height=7cm]{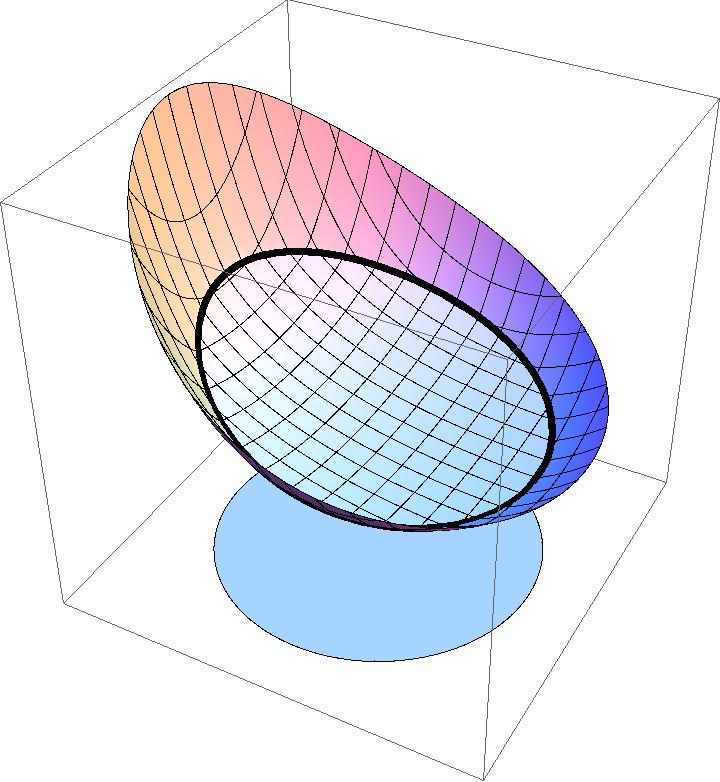}}
\caption{Modification of SK--metrics}
\end{figure}

{\bf Proof.}
By Theorem \ref{thm:existence} there exists a {\it unique} regular conformal metric $\nu(z)
\, |dz|$ with constant curvature $-1$ in $K$ such that $\nu$ is continuous on
$\overline K$ and $\nu(\xi)=\lambda(\xi)$ for all $\xi \in \partial K$. 
Exercise \ref{sec:hyperbolic}.\ref{exe:ex} shows $\lambda(z) \le \nu(z)$ for all $z \in K$.
Thus, by the Gluing Lemma, 
\begin{equation*}
M_K\lambda(z):= \begin{cases} \nu(z)  \quad  & \text{for } z \in K  \\
\lambda(z)  \quad& \text{for } z \in G \backslash K
\end{cases}
\end{equation*}
induces a (unique) SK--metric $M_K\lambda(z) \, |dz|$ with the desired properties.
\hfill{$\blacksquare$}

\begin{remark} \label{rem:mod}
Note that $M_K[\lambda] \ge \lambda$ and $M_K[\lambda] \ge M_K[\mu]$ if
$\lambda \ge \mu$. 
\end{remark}

\begin{definition}[Perron family]
A  family $\Phi$ of (densities of) SK--metrics on $G$ is called a Perron family, if the following conditions are satisfied: 
\begin{itemize}
\item[(i)]  If $\lambda \in \Phi$ and $\mu \in \Phi$, then $\sigma \in \Phi$, where $\sigma(z) := \max \{\lambda(z), \mu(z)\}$. 
\item[(ii)] If $\lambda  \in \Phi$, then $M_K\lambda \in
  \Phi$ for any open disk $K$ compactly contained in~$G$.
\end{itemize}
\end{definition}

\begin{theorem} \label{thm:perron0}
Let $\Phi$ be a Perron family  of SK--metrics  on $G$.
If $\Phi \not= \emptyset$, then 
 $$\lambda_{\Phi}(z):=\sup_{\lambda \in \Phi} \lambda(z)$$ induces  a regular conformal metric  $\lambda_{\Phi}(z)\, |dz|$  of constant curvature $-1$ on $G$.  
 \end{theorem}
 
{\bf Proof.}
 Fix $z_0 \in G$ and choose an open disk $K_R(z_0)$ in
$G$. Then $(T^* \lambda)(0)\le \lambda_{\D}(0)=2 $ for $T(z)=R z+z_0$
by the Fundamental Theorem, so $\lambda(z_0) \le 2/R$ for every $\lambda \in
 \Phi$ and $\lambda_{\Phi}$ is well--defined.

\smallskip

It suffices to show that $\lambda_{\Phi}(z) \, |dz|$ is regular and has constant
curvature $-1$ in every open disk  which is compactly contained in $G$.
 Now fix such an open disk $K$ and pick $z_0 \in K$.
Let $(\lambda_n) \subset \Phi$ such that
$\lambda_n(z_0) \to \lambda_{\Phi}(z_0)$ and let $\tilde{\lambda}_n \in \Phi$ denote the
modification of the SK--metric  $\max\{  \lambda_1, \ldots, \lambda_n\} \in \Phi$ on $K$.
Then $\tilde{\lambda}_n(z_0) \to \lambda_{\Phi}(z_0)$
and $\tilde{\lambda}_n(z) \, |dz|$ has curvature $-1$ in $K$. Since
the sequence $(\tilde{\lambda}_n)$ is monotonically increasing, it converges to a regular conformal metric
$\tilde{\lambda}(z) \, |dz|$ with curvature $-1$ on $K$ by Lemma \ref{lem:monotone}.

\smallskip

We claim that $\tilde{\lambda}=\lambda_{\Phi}$ in $K$.
By construction, $\tilde{\lambda} \le \lambda_{\Phi}$ in $K$
and $\tilde{\lambda}(z_0)=\lambda_{\Phi}(z_0)$. Assume, for a contradiction, that
$\tilde{\lambda}(z_1) <\lambda_{\Phi}(z_1)$ for some $z_1 \in K$. Then 
$\tilde{\lambda}(z_1)<\lambda(z_1)$ for some $\lambda \in \Phi$. Let $\mu_n
\in \Phi$ denote the modification of $\max\{\lambda,\tilde{\lambda}_n\}\in \Phi$ in
$K$. Then $\mu_n(z_0) \to \lambda_{\Phi}(z_0)$, $\mu_n \ge \tilde{\lambda}_n$ in $K$
and $(\mu_n)$ is monotonically increasing. By Lemma \ref{lem:monotone} the sequence
$(\mu_n)$ converges to a regular conformal metric $\mu(z) \, |dz|$ with
curvature $-1$ on $K$ with $\tilde{\lambda}(z_0)=\mu(z_0)$ and
$\tilde{\lambda} \le \mu$ in $K$. Lemma \ref{lem:equality} shows that
$\tilde{\lambda}=\mu$ in $K$, which contradicts $\tilde{\lambda}(z_1)<\lambda(z_1)
\le \mu_n(z_1) \to \mu(z_1)$.~\hfill{$\blacksquare$}

\subsection{The hyperbolic metric: Definition and basic properties}

The family of all SK--metrics on a domain $G$ is clearly a Perron family. Thus
we obtain as a special case of Theorem \ref{thm:perron0} the following result.

\begin{theorem} \label{thm:perron}
Let $\Phi_G$ be the Perron family  of all SK--metrics  on $G$.
If $\Phi_G \not= \emptyset$, then 
 $$\lambda_G(z):=\sup_{\lambda \in \Phi_G} \lambda(z)$$ induces  a regular conformal metric  $\lambda_{G}(z)\, |dz|$  of constant curvature $-1$ on $G$.  In particular, $\lambda_G(z)\, |dz|$ is the (unique) maximal SK--metric on $G$.
\end{theorem}

\begin{remark} 
 We call  $ \lambda_G(z)\, |dz|$ the hyperbolic metric on $G$. \label{def:hypi}
By Theorem \ref{thm:hyp} every domain with at least two boundary points
carries a hyperbolic metric. Note the obvious, but important monotonicity property
$\lambda_{D} \le \lambda_{G}$ if $G \subseteq D$.
\end{remark}

The Fundamental Theorem shows that the hyperbolic metric of the unit disk is
given by
$$ \lambda_{\D}(z) \, |dz| =\frac{2\, |dz|}{1-|z|^2} \, ;$$
 the hyperbolic metric on the punctured unit disk $\D'$
is 
$$ \lambda_{\D'}(z) \, |dz|=\frac{|dz|}{|z| \log (1/|z|)} \, , $$
see Exercise \ref{sec:ahlfors}.\ref{exe:punc} and Theorem \ref{thm:perron}.

\medskip

In general, however, it is very difficult to find an explicit formula for the
hyperbolic metric of a given domain. It is therefore important to obtain good {\it estimates}
for the hyperbolic metric. We shall use the Gluing Lemma for this purpose and
 start with a lower bound for the hyperbolic metric of the twice--punctured
plane $\C''$.

\begin{theorem} \label{thm:min}
Let
$$ \log R:=\frac{1}{\min \limits_{|z|=1} \lambda_{\C''}(z)} \, .$$
Then
$$ \lambda_{\C''}(z) \ge \frac{1}{|z| \left( \log R+\big|\log|z|\big| \right) } \, , \qquad z \in \C'' \, .$$
\end{theorem}

We shall later find the exact value of $R$, see Exercise \ref{sec:hyperbolic}.\ref{exe:value} and Corollary \ref{cor:val}.

\medskip

{\bf Proof.}
We first consider the case $z \in \D'$. By definition of $R$,
$$ \lambda_{\C''}(z) \ge  \frac{1}{\log R}=\frac{1}{|z| \log(R/|z|)}=\lambda_{\D_R'}(z) \quad \text{ for } z \in \partial \D  \, .
$$
Thus the Gluing Lemma guarantees that
$$\sigma(z):=\begin{cases} \,   \max \{\lambda_{\C''}(z), \lambda_{\D_R'}(z)\}   & \hspace{3mm} \, \text{for }  z \in \D' \, , \\[2mm]
                       \,        \lambda_{\C''}(z)         & \hspace{3mm} \,
                       \text{for } z \in \C''\backslash \D',  
          \end{cases}
$$
induces an SK--metric on $\C''$, so $\sigma(z) \le \lambda_{\C''}(z)$, which implies 
$$\lambda_{\C''}(z) \ge \lambda_{\D_R'}(z)=\frac{1}{|z| \log(R/|z|)} \quad
\text{ for }   z \in \overline{\D} \backslash \{ 0,1\}\, . $$
In a similar way (see Proposition \ref{prop:1}) we get 
$$\lambda_{\C''}(z) \ge \lambda_{\Delta_{1/R}}(z)
=\frac{1}{|z| \log(R \, |z|)} \quad \text{ for } z \in \C'' \backslash \D' \, . $$
Combining both estimates completes the proof.
\hfill{$\blacksquare$}

\medskip

It is now easy to show that
the behavior of the hyperbolic metric $\lambda_G(z) \, |dz|$ 
near an isolated boundary point 
mimics the behavior of the hyperbolic metric
$\lambda_{\D'}(z) \, |dz|$ of the punctured unit disk near the origin. 

\begin{corollary} \label{cor:is}
Let $G$ be a domain with an isolated boundary point $z_0 \in \C$ and hyperbolic metric $\lambda_{G}(z) \, |dz|$. Then 
$$ \lim \limits_{z \to z_0} \left( |z-z_0| \log \frac{1}{|z-z_0|} \right) \lambda_G(z)=1 \, . $$
\end{corollary}

{\bf Proof.}
We may assume $z_0=0$, $\D\backslash \{ 0\}  \subset G \backslash \{0 \}$ and $1 \in \partial G$. Then $\lambda_{\C''}(z) \le \lambda_{G}(z) \le \lambda_{\D'}(z)$ for any $z \in \D'$, 
so Theorem \ref{thm:min} leads to  
 $$ \frac{1}{|z| \log(R/|z|)} \le \lambda_{\C''}(z) \le 
\lambda_{G}(z) \le \lambda_{\D'}(z)=\frac{1}{|z| \log (1/|z|)} \, , \qquad z \in \D' \,,
$$
which proves the Corollary. \hfill{$\blacksquare$}

\subsection{Completeness}

Given a  conformal metric  $\lambda(z)\, |dz|$  on $G$ we define an associated
distance function by
\begin{equation*}
\d(z_0,z_1):= \inf_{\gamma} L_{\lambda}(\gamma) =\inf_{\gamma} \int_{\gamma} \lambda (z) \, |dz|\, ,
\end{equation*} \label{def:dist}
where the infimum is taken over all paths $\gamma$ in $G$ joining $z_0$ and $z_1$.

\begin{remark}
If $\lambda(z)\, |dz|$ is a conformal metric on $G$, then
 $(G, \d)$  is a metric space, see Exercise \ref{sec:hyperbolic}.\ref{exe:met}. 
\end{remark}

For instance, if we take the euclidean metric $\lambda(z) \, |dz|:=|dz|$
 on $G=\C$, then $\d$ is the euclidean distance, $\d(z_0,z_1)=|z_0-z_1|$.
As another example, we briefly discuss the distance induced by the hyperbolic metric $\lambda_{\D}(z) \, |dz|$ of the unit disk.

\begin{example}\label{ex:hypo}
We determine the distance $\d(z_0,z_1)$  associated to the hyperbolic metric $\lambda_{\D}(z) \, |dz|$ on the unit disk $\D$.
We first deal with the special case $z_0=0$ and $z_1 \in \D \backslash \{ 0 \}$.
Let $\gamma : [0,1] \to \D$ be a path connecting $0$ and $z_1$ with
$\gamma(t)\not=0$ for all $t \in (0,1]$.
 Then 
\begin{eqnarray*}
 L_{\lambda_{\D}}(\gamma) &=&
\int \limits_{0}^{1} \frac{2\, |\gamma'(t)|}{1-|\gamma(t)|^2} \, dt
\ge
 \int \limits_{0}^1 \frac{2\, \frac{d}{dt} |\gamma(t)|}{1-|\gamma(t)|^2} \, dt
\overset{s=|\gamma(t)|}{=}
\int \limits_{0}^{|z_1|} \frac{2\, ds}{1-s^2} = \log
\left( \frac{1+|z_1|}{1-|z_1|} \right) \, ,
\end{eqnarray*}
with equality if and only if $\gamma(t) \equiv z_1\,  t$. Thus
$$\text{\rm d}_{\lambda_{\D}}(0,z_1)= \log \left( \frac{1+|z_1|}{1-|z_1|} \right)\, .$$

\medskip

Now, let $z_0$ and $z_1$ be two distinct points in $\D$
and $T$ a conformal self--map  of $\D$ with  $T(z_0)=0$. By Exercise \ref{sec:curvature}.\ref{exe:1} $ \lambda_{\D}(T(z)) \, |T'(z)|=\lambda_{\D}(z)$,
so (\ref{eq:length}) implies $L_{\lambda_{\D}}(\gamma)=L_{\lambda_{\D}}(T \circ \gamma)$ for each path $\gamma \subset \D$. This leads to  
$$ \text{\rm d}_{\lambda_{\D}}(z_0,z_1)=\text{\rm d}_{\lambda_{\D}}(T(z_0),T(z_1))=\text{\rm d}_{\lambda_{\D}}\left(0, \frac{z_1-z_0}{1-\overline{z_0} \, z_1} \right)
= \log \left( \frac{1+\left| \displaystyle \frac{z_1-z_0}{1-\overline{z_0} \, 
      z_1}\right|}{1-\displaystyle \left|
    \frac{z_1-z_0}{1-\overline{z_0} \, z_1} \right|} \right)\,  .$$ 
\end{example}

If $\lambda(z) \, |dz|$ is a conformal metric on $G$, we can  equip
$G$ with the ordinary euclidean metric and get the metric space $(G,|\cdot|)$ and we also can equip $G$ with the distance $\d$ coming from $\lambda(z) \, |dz|$ and get the metric space $(G,\d)$. Topologically, these two metric spaces are equivalent.

\begin{proposition} \label{prop:topology}
Let $\lambda(z) \, |dz|$ be a conformal metric on $G$. Then the two metric spaces $(G,| \cdot |)$ and $(G,\d)$ have the same open and compact sets.
\end{proposition}

{\bf Proof.}
We will prove that the identity maps (i) ${\rm Id}:(G, | \cdot|) \to (G, \d)$
and (ii) ${\rm Id}:(G,
\d) \to  (G, | \cdot|)$ are continuous.

\smallskip 

(i) Pick $z_0 \in G$ and choose $\varepsilon>0$. Since $\lambda$ is continuous
on $G$, there exists a $\delta>0$ such that $K:=K_{\delta}(z_0)$ is compactly
contained in $G$ and
$$\delta \cdot \max \limits_{z \in K}\lambda(z) < \varepsilon\, .$$
Thus $\d(z, z_0) \le \int_{\gamma} \lambda(z) \, |dz| \le \delta \max_{z \in K} \lambda(z)< \varepsilon$ for all $z \in K$, where $\gamma$ is 
the straight line connecting $z_0$ and $z$. 

\smallskip

(ii) Fix $z_0 \in G$ and let $(z_n)$ be a sequence in $G$ with $\d(z_n,z_0) \to
0$ as $n \to \infty$. Further, let $\varepsilon >0$ such that
$K:=\overline{K_{\varepsilon}(z_0)} \subset G$ and set $$m:= \min_{z \in K} \lambda(z)\, .$$ 
Then there is an integer $N$ such that 
$$\d(z_n, z_0) < m\, \varepsilon$$
for all
$n\ge N$. We now assume that we can find some $j  \ge N$ such that $| z_j -z_0| >
\varepsilon$. In particular,  $z_j
\not \in K$. Let $\gamma:[0,1] \to G $ be a path joining $z_0$ and $z_j$. We
set $\tau:=\inf\,  \{ \, t \in [0,1]\, : \, \gamma(t) \not \in K\, \}$. Then
$\tilde{\gamma}:=\gamma_{|_{[0,\tau]}}$ is a path in $K$ and we conclude
$$ L_{\lambda}(\gamma)=\int_{\gamma} \lambda(z)\, |dz| \ge m \int_{\tilde{\gamma}}
|dz| \ge m\, \varepsilon\,.  $$ 
This shows $\d(z_j, z_0) \ge m \, \varepsilon$, contradicting our
hypothesis. Thus $| z_n -z_0| \to 0$  whenever $\d(z_n,z_0) \to
0$. \hfill{$\blacksquare$}

\medskip

However, even though 
the two metric spaces $(G,|\cdot|)$ and $(G,\d)$ are topologically equivalent, they are 
in general not {\it metrically} equivalent. For instance, in  $(\D,|\cdot|)$ the boundary $\partial \D$ has distance $1$ from the center $z=0$, but in $(\D,\text{d}_{\lambda_{\D}})$ the boundary is infinitely away from $z=0$. This latter property turns out to be particularly useful, so we make the following definition.

\begin{definition}
A conformal metric  $\lambda(z)\, |dz|$ on $G$ is called complete (for $G$), if 
\begin{equation*}
\lim_{n \to \infty} \d (z_n,z_0)= +\infty
\end{equation*} 
for every sequence $(z_n) \subset G$ which leaves every compact subset of $G$ and for some (and therefore for every) point $z_0 \in G$. 
\end{definition}

\begin{remarks}\label{rem:complete}
\begin{itemize}
\item[(a)] If $\lambda(z) \, |dz|$ is a complete conformal metric for $G$, then
$(G,\d)$ is a complete metric space in the usual sense, i.e., every
Cauchy sequence in $(G,\d)$  converges in $(G,\d)$ and,
by Proposition \ref{prop:topology}, then also in $(\C,| \cdot|)$ to some point {\it in} $G$.
\item[(b)] If $\lambda(z) \, |dz|$ is a complete conformal metric on $G$ and $\mu(z) \ge
\lambda(z)$ in $G$, then $\mu(z) \, |dz|$ is also complete (for $G$).
\end{itemize}
\end{remarks}

Obviously, the euclidean metric is complete for the complex plane $\C$ and 
 the hyperbolic metric $\lambda_{\D}(z) \, |dz|$ is complete for the unit disk $\D$, see Example \ref{ex:hypo}. 

\begin{example}\label{ex:complete1}
The hyperbolic metric $\lambda_{\D'}(z) \, |dz|$ on the punctured disk  $\D'$ is
  complete. By Exercise \ref{sec:curvature}.\ref{exe:2}  we have 
$$\lambda_{\D}(z)\, |dz|= \left( \pi^* \lambda_{\D'}\right)(z) \, |dz|\, ,$$
where $\pi : \D \to \D'$, $z \mapsto \text{ \rm exp}(- (1+z)/(1-z) )$, is  locally one-to-one and onto. Now pick $z_0, z_1 \in \D'$ and let $\gamma : [a,b] 
\to\D'$ be a path joining $z_0$ and $z_1$. Fix a branch $\pi ^{-1}$ of the inverse function of $\pi$ defined originally only in a neighborhood of $z_0=\gamma(a)$. Then $\pi^{-1}$ can clearly be analytically continued along $\gamma$ and we get the path  $\tilde{\gamma}:= \pi^{-1} \circ \gamma \subset \D$. Thus
we obtain
$$\text{\rm d}_{\lambda_{\D}}( \pi^{-1}(z_0), \pi^{-1}(z_1)) \le
L_{\lambda_{\D}}( \tilde{\gamma})= L_{\pi^*\lambda_{\D'}}(\tilde{\gamma}) \stackrel{\text{(\ref{eq:length})}}{=}
  L_{\lambda_{\D'}}(\gamma)\, , $$  
so, in particular,
$$ \text{\rm d}_{\lambda_{\D}}( \pi^{-1}(z_0), \pi^{-1}(z_1)) \le \text{\rm
  d}_{\lambda_{\D'}}(z_0, z_1)\, .$$
The desired result follows from the latter inequality
since $|\pi^{-1}(z_1)| \to 1$ if $z_1 \to 0$ or $|z_1| \to 1$.
\end{example}

\begin{theorem} \label{thm:hyp20}
The hyperbolic metric $\lambda_{\C''}(z) \, |dz|$ is complete. In particular,
the hyperbolic metric $\lambda_G(z)\, |dz|$  for any  domain $G$ with at least
two boundary points is complete.
\end{theorem}

{\bf Proof.} By Proposition \ref{prop:topology} 
  it suffices to  show that for
fixed $z_0 \in \C''$ we have $ \text{d}_{\C''}(z_n,z_0) \to
+\infty$ for any sequence $(z_n) \subset \C''$ which converges either to $0$, $1$ or
$\infty$.

\smallskip

If $z_n \to 0$, then by making use of  the estimate
$$ \lambda_{\C''}(z)  \ge \frac{1}{|z| \log(R/|z|)}=\lambda_{\D_{R}'}(z) \, ,\quad
 z \in \D' \, , $$ for some $R>0$ (see Theorem \ref{thm:min}),  we deduce that $\text{d}_{\C''}(z_n,z_0) \to +\infty$, since $\lambda_{\D'_R}(z) \, |dz|$ is complete for $\D'_R$, cf.~Example \ref{ex:complete1}.

\smallskip

In order to deal with the other cases, we note that  Exercise  \ref{sec:hyperbolic}.\ref{exe:wert1/2} gives us
$$ \lambda_{\C''}(z)= \lambda_{\C''}(T(z))\, |T'(z)| \, , \quad z \in \C''\, ,$$
for $T(z)=1/z$ and $T(z)=1-z$. Hence in view of (\ref{eq:length}),
\begin{eqnarray*}
\text{d}_{\C''}(z_1,z_0)&=&\text{d}_{\C''}\left(\frac{1}{z_1},\frac{1}{z_0}\right)\, ,  \\[2mm]
\text{d}_{\C''}(z_1,z_0)&=&\text{d}_{\C''}(1-z_1,1-z_0)\, . 
\end{eqnarray*}
This immediately shows $\text{d}_{\C''}(z_n,z_0) \to + \infty$ if $z_n \to
\infty$ and $\text{d}_{\C''}(z_n,z_0)\to + \infty$ if $z_n \to 1$. The
desired result follows.
\hfill{$\blacksquare$}

\medskip

As an application,
we now make use of the completeness of the hyperbolic metric to prove Montel's
extension of the theorems of Picard.

\begin{theorem}[The Big Montel Theorem]
The set of all holomorphic functions $f : \D \to \C''$ is a normal family.
\end{theorem}

{\bf Proof.} Let $\lambda:=\lambda_{\C''}$.
For every holomorphic function $f : \D \to \C''$, the Fundamental Theorem shows
 $f^*\lambda \le \lambda_{\D}$, so
$$ \text{d}_{\lambda} (f(z_1),f(z_2))\le \text{d}_{\lambda_{\D}}(z_1,z_2) \, , \qquad z_1, z_2
\in \D \, . $$
Thus the family ${\cal F}:=\{f : (\D,\text{d}_{\lambda_{\D}}) \to (\C'',\d) \text{ holomorphic}\}$
is equicontinuous at each point of $\D$. 
We distinguish two cases.

\smallskip

1.~Case: For each $z \in \D$ the set $\{f(z) \, : \, f \in {\cal F}\}$
has compact closure in $(\C'',\d)$. Then Ascoli's theorem (see \cite[Thm.~0.4.11]{Edw}) shows that every
sequence in ${\cal F}$ has a subsequence which converges locally uniformly in
$\D$ to some function in ${\cal F}$.

\smallskip

2.~Case: There exists a point $z_0 \in \D$ such that the closure of
$\{f(z_0) \, : \, f \in {\cal F}\}$ is not compact in $(\C'',\d)$.
Thus there is a sequence $(f_j) \subset {\cal F}$ such that
 $(f_j(z_0))$ converges to either $0$ or $1$ or $\infty$.
It suffices to consider the case $f_j(z_0) \to 0$.
Since $\d(f_j(z),f_j(z_0)) \le \text{d}_{\lambda_{\D}}(z,z_0)$ and
$\lambda_{\C''}(z) \, |dz|$ is complete, it follows that $f_j(z) \to 0$ for each $z \in \D$.
This convergence is actually locally uniform in $\D$, since otherwise one could
find  points $z_j \in \D$ such that $z_j \to z_1 \in \D$ and $|f_j(z_j)| \ge \eps>0$.
But then, once again, the completeness of $\lambda_{\C''}(z) \, |dz|$ would imply $+\infty
\leftarrow\d(f_j(z_j),f_j(z_0)) \le \text{d}_{\lambda_{\D}}(z_j,z_0) \to
\text{d}_{\lambda_{\D}}(z_1,z_0)<+\infty$. Contradiction!
\hfill{$\blacksquare$}

\section*{Exercises for Section 4}

\begin{enumerate}
\item \label{exe:4.1}  Show that there exists no SK--metric on the punctured plane $\C'$.

\item \label{exe:ex}
Let $\lambda(z) \, |dz|$ be a conformal pseudo--metric on $G$. Show that the following are equivalent:
\begin{itemize}
\item[(a)] $\lambda(z) \, |dz|$ is an SK--metric on $G$.
\item[(b)] If $D$ is  compactly contained in $G$ and $\mu(z) \, |dz|$ is a   regular conformal metric with constant curvature $-1$ on $D$ such that
$$ \limsup \limits_{z \to \xi} \frac{\lambda(z)}{\mu(z)} \le 1 \quad \text{ for every } \xi \in \partial D \, , $$
then $\lambda \le \mu$ in $D$.
\end{itemize}

\item \label{exe:4.3} Show that the pullback of an SK--metric under a
  non--constant analytic map is again an SK--metric.

\item \label{exe:4.4a} Let $\lambda(z) \, |dz|$ be an SK--metric on $\D'$ such
  that
$$ \limsup_{z \to 0} \lambda(z) \, |z|^{\alpha} < + \infty$$
for some $\alpha <1$. Let 
$$\lambda_{\alpha}(z)= \frac{2\, (1-\alpha) |z|^{- \alpha}}{1- |z|^{2\, (1-
    \alpha)}}\, .$$
Show that 
$$  \lambda(z) \le \lambda_{\alpha}(z)$$
for all $z \in \D'$. 

\item \label{exe:4.4} Let  
$$\lambda_{n}(z) =\frac{2\, \left(1+\frac{1}{n} \right)  \, |z|^{\frac{1}{n}}}{1- |z|^{2
  (1+\frac{1}{n})}} \, .$$
Show that $(\lambda_{n})$ is a monotonically increasing sequence of
densities of SK--metrics in $\D$, whose limit $\lambda$ does {\it not} induce an SK--metric on $\D$.

\item \label{exe:heins} Use the Gluing Lemma to prove the following theorem of M.~Heins
(\cite[Theorem 18.1]{Hei62}):
{\it Let $G$ be a domain with a hyperbolic metric $\lambda_{G}(z) \, |dz|$ and let $D \subseteq G$ denote a component of the complement of a compact subset 
 $K \subset G$. Then $\lambda_{G}/\lambda_{D}$ has a positive lower bound on
$D \backslash V$ for any neighborhood $V$ of $K$.}

\item \label{exe:wert1/2}
Let $T$ be one of the M\"obiustransformations
$$ 1/z \, , \qquad 1-z \, , \qquad \frac{z}{z-1} \, .$$
Show that $T^* \lambda_{\C''}=\lambda_{\C''}$ and deduce
$\lambda_{\C''}(1/2)=4 \, \lambda_{\C''}(-1)$.

\item \label{exe:schwarzpicard}

(The Schwarz--Picard Problem on the sphere)

A point $z_0 \in \hat{\C}:=\C \cup \{ \infty\}$ is called a conical singularity
of order $\alpha \le 1$ of an SK--metric $\lambda(z) \, |dz|$ defined in a punctured neighborhood of $z_0$ if 
$$\begin{array}{rcll}
\log \lambda(z) &\!\!=\!\!&  \begin{cases} -\alpha \log |z-z_0|+O(1) \hspace{3.2cm} & \text{ as } z \to z_0 \not=\infty \\ -(2-\alpha) \log|z|+O(1) & \text{ as } z  \to z_0=\infty \,  
\end{cases}\, & \text{ if } \alpha <1  \\[6mm]
\log \lambda(z)  & \!\!=\!\!& \begin{cases} - \log |z-z_0|-\log \left( -\log |z-z_0| \right)+O(1) & \text{ as } z \to z_0 \not=\infty \\ - \log|z|+ \log \log |z|+O(1) & \text{ as } z  \to z_0=\infty \, 
\end{cases}\, & \text{ if } \alpha=1 \, .
\end{array}$$
Now, let $n \ge 3$ distinct points $z_1, \ldots \, ,z_{n-1},\infty \in \hat{\C}$ 
and real numbers $\alpha_1, \ldots , \alpha_n \le 1$ be given such that
\begin{equation} \label{eq:gaussbonnet}
 \sum \limits_{j=1}^n \alpha_j >2 \, .
\end{equation}
Choose $\delta>0$ such that
$$ \min \limits_{k\not=j} |z_j-z_k|>\delta \quad \text{ and } \quad  |z_j|<1/\delta \, \quad \text{ for }  \quad j=1, \ldots, n-1 \, .$$
Let
\begin{align*}
&\lambda_j(z)= \begin{cases} 
\displaystyle 2\, (1- \alpha_j) \, 
\frac{\delta^{1-\alpha_j}\, |z-z_j|^{-\alpha_j}}{\delta^{2(1-\alpha_j)} -|z-z_j|^{2(1-\alpha_j)}} \qquad \qquad &
\text{if } \,  \alpha_j<1\\[5mm]
\displaystyle \frac{1}{|z-z_j|\, \log(\delta/|z-z_j|)} & \text{if } \, \alpha_j=1
    \end{cases}
\intertext{for $j=1, \ldots, n-1$, and}
&\lambda_n(z)= \begin{cases} 
\displaystyle 2\, (1- \alpha_n)\,  \frac{\delta^{1-\alpha_n}   |z|^{-
    \alpha_n}} {\delta^{2(1-\alpha_n)} \, |z|^{2(1-\alpha_n)}-1} \qquad \qquad &
\text{if } \,  \alpha_n<1\\[5mm]
\displaystyle \frac{1}{|z|\, \log(\delta \, |z|)} & \text{if } \, \alpha_n=1\, .
    \end{cases}
\end{align*}

\begin{itemize}
\item[(a)] Denote by $\Phi$ the family of all densities of SK--metrics
  $\lambda(z) \, |dz|$ on the $n$--punctured plane $G:=\C \backslash \{z_1, \ldots ,z_{n-1}\}$ such that
$$\limsup_{z \to z_j} \frac{\lambda(z)}{\lambda_j(z)} \le 1$$
for all $j=1, \ldots, n$. 
Show that $\Phi$ is a non--empty  Perron family.
\item[(b)] Show that $\lambda_{\Phi}(z):=\sup_{\lambda \in \Phi} \lambda(z)$ 
induces the uniquely determined regular conformal metric of constant curvature $-1$ on $G$
with conical singularities of order $\alpha_j$ at $z_j$ for $j=1, \ldots, n$.
\end{itemize}
(Note: Condition (\ref{eq:gaussbonnet}) is also necessary (Gauss--Bonnet).)

\item \label{exe:app}
Let $\Omega_n:=\C'' \backslash \{e^{2 \pi i k/n} \, : \, k=1, \ldots , n-1\}$
for each positive integer n and let $f_n(z)=z^n$. Show that
$\lambda_{\Omega_n}(z) \, |dz| = (f_n^*\lambda_{\C''})(z) \, |dz|$.


\item \label{exe:value}
The aim of this exercise is to show that $\min \limits_{|z|=1} \lambda_{\C''}(z)=\lambda_{\C''}(-1)$.

For $\eta \in \partial \D$ let $\lambda_{\eta}(z):=\lambda_{\C \backslash \{0, \eta\}}(z)$.
\begin{itemize}
\item[(a)] Show that $\lambda_{\overline{\eta}}(z)=\lambda_{\eta}(\overline{z})$.
\item[(b)] Use the Gluing Lemma to show that $\lambda_{\overline{\eta}}(z) \le 
\lambda_{\eta}(z)$ for all $\Im z>0$ if $\Im \eta>0$.

(Hint: Consider $\max \{ \lambda_{\overline{\eta}}(z),\lambda_{\eta}(z)\}$
on the upper half plane and $\lambda_{\eta}(z)$ in the lower half plane.)
\item[(c)] Now fix $\theta \in (-\pi,0)$ and let $\eta:=e^{-i \theta/2}$.
Verify that
$$ \lambda_{\C''}(-e^{i \theta})\overset{(a)}{=}\lambda_{\eta}(-\overline{\eta})
\overset{(b)}{\ge} \lambda_{\overline{\eta}}(-\overline{\eta})\overset{(a)}{=}\lambda_{\C''}(-1) \,.$$
\item[(d)] Show in addition 
that $\theta \mapsto \lambda_{\C''}(r e^{i \theta})$
is strictly decreasing on $(0,\pi)$ and strictly increasing on $(-\pi,0)$ for each $r>0$.
\end{itemize}
\item \label{exe:met}
Let $\lambda(z) \, |dz|$ be a conformal metric on $G$. Show that
$(G,\d)$ is a metric space.
 \end{enumerate}

\section*{Notes}

Most of the material of this section and much more can be found in M.~Heins
\cite{Hei62}, see also S.~Smith \cite{Smi86}. We again wish to emphasize the
striking similarities of regular conformal metrics of constant curvature $-1$ and SK--metrics
with  harmonic and subharmonic functions, see Ransford \cite{Ransford}.
Theorem \ref{thm:min} was proved by J.~A.~Hempel \cite{Hem79} (see also
D.~Minda \cite{Min87b}).
Different proofs for Corollary \ref{cor:is} were given by
 J.~Nitsche \cite{Nit57},  M.~Heins \cite[Section 18]{Hei62},  A.~Yamada \cite{Yam1}, S.~Yamashita \cite{Yam2} and D.~Minda \cite{Min97}. 
The metric space $(G,\d)$ is an example of a path metric space in the sense of Gromov \cite{Gro99}.
The converse of Remark \ref{rem:complete} (a) is the Hopf--Rinow Theorem. It says that  $(G,\d)$ is a complete metric space if and only if every closed and bounded set in $(G,\d)$ is compact. See Gromov \cite[Chapter 1]{Gro99} for a quick proof of the Hopf--Rinow Theorem. For more information about
normal families
of analytic functions and conformal metrics beyond Montel's Big Theorem we refer to 
Grauert and Reckziegel \cite{GR65}. The Schwarz--Picard Problem (Exercise \ref{sec:hyperbolic}.\ref{exe:schwarzpicard})  was first solved by Poincar\'e \cite{Poi1898} in the case $\alpha_1=\ldots=\alpha_n=1$.   The general case (even on a compact Riemann surface instead of the sphere) was treated for instance by Picard \cite{Pic1893,Pic1905}, 
Bieberbach \cite{Bie16}, Lichtenstein \cite{Li}, Heins \cite{Hei62}, McOwen \cite{McO88,McO93} and Troyanov \cite{Tro90}. 
See Lehto, Virtanen \& V\"ais\"ala \cite{Leh2},
Hempel \cite{Hem79}, Weitsman \cite{Wei79} 
and Minda \cite{Min87b} for Exercise
\ref{sec:hyperbolic}.\ref{exe:value}. There one can find much 
more information about monotonicity properties of the hyperbolic metric.

\section{Constant Curvature} \label{sec:constcurv}

The pullback $(f^*\lambda_{\D})(z) \, |dz|$ of the hyperbolic metric
$\lambda_{\D}(w) \, |dw|$ under a locally univalent analytic function $w=f(z)$
from $G$ to $\D$ is a
regular conformal metric $\lambda(z) \, |dz|$ on $G$ with constant curvature $-1$.
We shall prove the following local converse.

\begin{theorem}[Liouville]\label{thm:liouville}
Let $\lambda(z)\, |dz|$ be a regular  conformal metric of constant curvature $-1$ on a simply connected domain $G$. Then there exists a locally univalent analytic function $f: G \to \mathbb{D}$ such that
\begin{equation} \label{eq:liouville2} 
\lambda(z)  = \frac{2\, |f'(z)|}{1-|f(z)|^2} \, , \qquad z \in G \, .
\end{equation}
If $g:G  \to \mathbb{D}$ is another locally univalent analytic function, then 
\begin{equation*}  
\lambda(z)  = \frac{2\, |g'(z)|}{1-|g(z)|^2} \, , \qquad z \in G \, ,
\end{equation*}
 if and only if $g=T\circ f$, where $T$ is a conformal self--map of $\D$. 
 \end{theorem}

Thus on simply connected domains conformal metrics of constant curvature $-1$ 
and bounded locally univalent functions can be identified.
More precisely, if $G$ is a simply connected domain with $z_0 \in G$, let 
\begin{equation*}
\Lambda(G):=\{  \lambda  \, : \, \lambda(z) \, |dz|  \text{  regular conformal metric with  } \kappa_{\lambda} \equiv -1 \text{ on }  G \}
\end{equation*}
and
\begin{equation*}
{\cal H}_0(G):=\{  f\, : \, f:G \to \mathbb{D} \text{ analytic and locally
  univalent with } f(z_0)=0, f'(z_0)>0\} \, .
\end{equation*}  
Then the map 
$$\Psi: {\cal H}_0(G) \to \Lambda(G)\, , \qquad f \mapsto \lambda=\frac{2\ |f'|}{1-|f|^2}\, , $$
is one--to--one. 

\begin{remark} \label{rem:21}
The proof of Theorem \ref{thm:liouville} below will show that when  $G$ is not simply connected, then the 
function  $f$  can be analytically continued along
any path $\gamma \subset G$ and (\ref{eq:liouville2}) holds along the path $\gamma$.
We call the (multivalued and locally univalent) function $f$ the developing map of the constantly
curved metric $\lambda(z) \,
|dz|$. The developing map is uniquely determined by the metric up to
postcomposition with a conformal self--map of $\D$.
\end{remark}

\begin{example1}
The conformal metric
$$ \lambda(z) \, |dz|=\frac{|dz|}{\sqrt{|z|} \, (1-|z|)}$$
on the punctured unit disk $\D'$ has constant curvature $-1$.
Its developing map  (modulo normalization) is the multivalued analytic 
function $f(z)=\sqrt{z}$.
\end{example1}

\subsection{Proof of Theorem \ref{thm:liouville}}

The proof of Theorem \ref{thm:liouville} will be split into
several lemmas which are of independent interest. We start off with the
following preliminary observation, see Remark \ref{rem:representation}.

\begin{remark} \label{rem:34}
If $\lambda(z)\, |dz|$ is a regular conformal metric of constant curvature $-1$, then the function $\lambda$ is of class $C^{\infty}$. 
\end{remark}

In order to find for a
conformal metric $\lambda(z) \, |dz|$ with constant curvature $-1$ a
holomorphic function $f$ such that 
the representation formula (\ref{eq:liouville2}) holds, we first
associate to $\lambda(z) \, |dz|$ a {\it holomorphic} auxiliary function.

\begin{lemma} \label{lem:liou0}
Let $\lambda(z)\, |dz| $ be a regular conformal metric of constant curvature
$-1$  on $G$ and $u(z):=\log \lambda(z)$. Then  
$$
 \frac{\partial u^2}{\partial z^2}(z) - \left( \frac{\partial u}{\partial z}(z) \right)^2
$$
 is holomorphic in $G$.
\end{lemma}

{\bf Proof.}
Since $u$ is a solution to $\Delta u=e^{2 u}$, we can write
\begin{equation*} 
\frac{\partial^2 u}{\partial z \partial \overline{z}}(z)=\frac{1}{4} \, e^{2 u(z)}
\, .
\end{equation*}
In view of Remark \ref{rem:34}, we are allowed to
differentiate this identity with respect to $z$. Hence
$$ \frac{\partial^3 u}{\partial z^2 \partial \overline{z}}=\frac{1}{2} \,
e^{2 u} \, \frac{\partial u}{\partial z}=2 \, 
\frac{\partial^2 u}{\partial z \partial \overline{z}} \, \frac{\partial
u}{\partial z}=\frac{\partial}{\partial \overline{z}} \left[ \left( \frac{\partial
  u}{\partial z} \right)^2 \right] \, ,
$$
and therefore
$$ 
\frac{\partial}{\partial \overline{z}} \left[ \frac{\partial^2u}{\partial z^2}(z)-\left( \frac{\partial u}{\partial z}(z)
\right)^2 \right] \equiv 0 \, .$$
\hfill{$\blacksquare$}

\begin{definition}
Let $\lambda(z)\, |dz| $ be a regular conformal metric and $u(z):=\log
\lambda(z)$. Then the function
\begin{equation} \label{eq:schwarz}
S_{\lambda}(z):= 2 \left[ \frac{\partial u^2}{\partial z^2}(z) - \left( \frac{\partial u}{\partial z}(z) \right)^2 \right]
\end{equation}
is called the Schwarzian derivative of  $\lambda(z)\, |dz| $.
\end{definition}

The appearance of the constant $2$ in this definition is motivated by

\begin{example} \label{l1}
If $\lambda(z) \, |dz|$ has the form
$$ \lambda(z)=\frac{2 \, |f'(z)|}{1-|f(z)|^2} \, ,$$
 then the Schwarzian derivative of the metric $\lambda(z) \, |dz|$ is equal to 
 the Schwarzian derivative of the analytic function $f$, that is
\begin{equation} \label{eq:liouvivp0}
\displaystyle S_f(z):=\left(  \frac{f''(z)}{f'(z)} \right)'-\frac{1}{2}\left(  \frac{f''(z)}{f'(z)} \right)^2= 
 S_{\lambda}(z) \, . 
\end{equation}
\end{example}

Thus, to find for a conformal metric $\lambda(z) \, |dz|$ with
constant curvature $-1$ a holomorphic function $f$ such that (\ref{eq:liouville2})
holds, first one computes the Schwarzian derivative $S_{\lambda}$ and then solves
the Schwarzian differential equation (\ref{eq:liouvivp0}) for $f$. In order to
guarantee that this program works one needs to know that 
conformal metrics with constant curvature $-1$ are (modulo a normalization)
uniquely determined by their Schwarzian derivatives.

\begin{lemma}\label{lem:liou1}
Let $\lambda(z) \, |dz|$ and $\mu(z) \, |dz|$ be regular conformal metrics
of constant curvature $-1$ on $D$ such that
\begin{equation*}
S_{\lambda}(z)=S_{\mu}(z)
\end{equation*}
for every $z \in D$ and
\begin{equation*}
\lambda(z_0)=\mu(z_0) \quad \text{and} \quad \frac{\partial \lambda}{\partial z} (z_0) = \frac{\partial \mu}{\partial z} (z_0)
\end{equation*}
for some point $z_0 \in D$. Then $\lambda(z) = \mu(z)$ for every $z \in D$. 
\end{lemma}

{\bf Proof.} 
(a) \, Let $u(z)=\log \lambda(z)$ and $v(z)=\log \mu(z)$. Then
\begin{equation}\label{eq:nitsche1}
\frac{\partial u^2}{\partial z^2}(z) - \left( \frac{\partial u}{\partial z}(z) \right)^2=\frac{\partial v^2}{\partial z^2}(z) - \left( \frac{\partial v}{\partial z}(z) \right)^2
\end{equation}
for every $z \in D$ and
\begin{equation}\label{eq:nitsche2}
u(z_0)=v(z_0) \quad \text{and} \quad \frac{\partial u}{\partial z} (z_0) =
\frac{\partial v}{\partial z} (z_0) \, .
\end{equation}

 We first show that 
\begin{equation} \label{eq:nitsche3}
\frac{\partial^{k+j} u}{\partial z^k \partial \bar{z}^j}(z_0)=
\frac{\partial^{j+k} v}{\partial z^k \partial \bar{z}^j} v(z_0)
\end{equation}
for $k,j=0,1, \ldots$.
By Lemma \ref{lem:liou0} both sides of equation (\ref{eq:nitsche1}) are holomorphic
functions in $D$ and it follows inductively from (\ref{eq:nitsche1}) and
(\ref{eq:nitsche2}) 
that (\ref{eq:nitsche3}) holds for $j=0$ and each $k=0,1, \ldots$. 
Since $u$ is real--valued  (\ref{eq:nitsche3}) is  valid  for $k=0$ and
every $j=0,1, \ldots$. 
 But then the partial differential equation
$\Delta u=e^{ 2u}$ implies that (\ref{eq:nitsche3}) holds for all $k,j=0,1, \ldots$.

\medskip

(b) \,
 We note that $ v_1(z):=e^{-u(z)}$ and $v_2(z):=e^{-v(z)}$ are (formal) solutions
of the linear equation
$$ v_{zz}+\frac{S_{\lambda}(z)}{2} \, v=0 \, .$$
Thus we consider the
auxiliary function (Wronskian)
$$ W(z):=\frac{\partial v_1}{\partial z}(z) \, v_2(z) -v_1(z) \, \frac{\partial
  v_2}{\partial z}(z) \, . $$
Then a straightforward calculation, using $S_{\lambda}=S_{\mu}$, shows 
$$\frac{\partial W}{\partial z}(z) =0 \, , $$
so that $W$ is antiholomorphic in $D$. But (a) implies that
$$ \frac{\partial^k W}{\partial \bar{z}^k}(z_0)=0$$
for each $k=0,1, \ldots$, so $W \equiv 0$ in $D$. Consequently,
$$ \frac{\partial}{\partial z} \left( \frac{v_1(z)}{v_2(z)} \right)
=\frac{W(z)}{v_2(z)^2} \equiv 0 \, ,$$
i.e.~$v_1/v_2$ is antiholomorphic in $D$. Since $v_1/v_2$ is real--valued
in $D$ it has to be constant $v_1(z_0)/v_2(z_0)=1$ by (\ref{eq:nitsche2}).
Therefore $u \equiv v$.
\hfill{$\blacksquare$}

\medskip

{\bf Proof of Theorem \ref{thm:liouville}.}
 Let $u(z):=\log \lambda(z)$. Lemma
\ref{lem:liou0} implies that $S_{\lambda}(z)$ is a holomorphic function in $G$.
Fix $z_0 \in G$.
Then, as it is well--known (see \cite[p.~53]{Leh} or \cite{Lai93}), the initial value problem 
\begin{equation} \label{eq:liouvivp}
\begin{array}{rl}
& \displaystyle S_f(z) = 
 S_{\lambda}(z) \, \\[4mm]
& f(z_0)=0 \, , \quad  \displaystyle f'(z_0)=\frac{\lambda(z_0)}{2} \, , \displaystyle \quad f''(z_0)=\frac{\partial
  \lambda}{\partial z}(z_0) \, ,
\end{array}
\end{equation}
has a unique meromorphic solution $f : G
\to \C$. Let $D$ be the connected component
of the set $\{z \in G \, : \, f(z) \in \D\ \text{ and } f'(z) \not=0\}$ which contains $z_0$ and 
let 
$$\mu(z):= \frac{2\, |f'(z)|}{1-\displaystyle |f(z)|^2}\, , \qquad z \in D
\, . $$
Example \ref{l1} and the choice of the initial conditions for $f$ 
guarantee  that the  hypotheses of Lemma \ref{lem:liou1}
are satisfied in $D$, so $\lambda(z)=\mu(z)$
for each $z \in D$. It is not difficult to show that this implies
that $D$ is closed (and open) in $G$, i.e., $D=G$ and
(\ref{eq:liouville2}) holds for each $z \in G$. 

\smallskip

Now assume that 
$$ \lambda(z)=\frac{2\, |g'(z)|}{1- |g(z)|^2} \, , \qquad z \in G \,
, $$
for some locally univalent holomorphic  function  $g : G \to \D$. There is a conformal self--map $T$ of $\D$ with $T(g(z_0))=0$ and $(T \circ g)'(z_0)>0$.
Then, by Exercise \ref{sec:curvature}.\ref{exe:1},
$$ \frac{2\, |(T \circ g)'(z)|}{1-\, |(T \circ g)(z)|^2} 
=\frac{2\, |T'(g(z))|}{1-\, |T (g(z))|^2} \,  |g'(z)|= 
\frac{2\, |g'(z)|}{1- |g(z)|^2}=\lambda(z)
\, ,
\qquad z \in G \, .
 $$
Thus  $T \circ g$ is a solution of the initial value problem (\ref{eq:liouvivp})
and therefore has to be identical to $f$ by uniqueness of this solution, so
$f= T \circ g$. \hfill{$\blacksquare$}

\subsection{Applications to the hyperbolic metric}

Theorem \ref{thm:liouville} is of fundamental importance and has many
applications. For instance, it can be used to analyze the Schwarzian of the
hyperbolic metric near isolated boundary points.

\begin{theorem} \label{thm:sch}
Let $G$ be a domain with an isolated boundary point $z_0$ and hyperbolic
metric $\lambda_G(z) \, |dz|$. Then $S_{\lambda_G}$ has a pole of order $2$ at $z=z_0$ and
$$ S_{\lambda_G}(z)=\frac{1}{2 \, (z-z_0)^2}+\frac{C}{z-z_0}+ \ldots \, $$
for some $C \in \C$.
\end{theorem}

{\bf Proof.} 
We may assume $z_0=0$ and $\D' \subseteq G$.
Let $f_n(z)=z^n$ and $\lambda_n(z) \, |dz|:=\left(f^*_n\lambda_G\right)(z) \, |dz|$. Then
by Corollary \ref{cor:is},
$$\lambda_n(z)=\lambda_G(z^n)\, n\, |z|^{n-1}= \frac{1}{|z| \, \log(1/|z|)} \,
\mu(z^n)\, ,$$
where $\mu$ is a continuous function at $z=0$ with $\mu(0)=1$.
This shows that $\lambda_n \to \lambda_{\D'}$ locally uniformly in $\D'$ as $n
\to \infty$.
Exercise \ref{sec:constcurv}.\ref{exe:schwarztransform} and
Liouville's Theorem imply
$$ S_{\lambda_G}(z^n) \, n^2 z^{2n-2}+\frac{1-n^2}{2 \, z^2}= S_{\lambda_n}(z)
 \to S_{\lambda_{\D'}}(z)=\frac{1}{2 \, z^2}$$
locally uniformly in $\D'$. Now a comparison
of the  Laurent coefficients on both sides immediately
completes the proof. \hfill{$\blacksquare$}

\medskip

We now apply Liouville's Theorem to {\it complete} conformal metrics of
curvature $-1$.

\begin{theorem} \label{thm:liouvillecomplete}
Let $\lambda(z)\, |dz|$ be a regular  conformal metric of constant curvature
$-1$ on $G$. Then the following are equivalent.
\begin{itemize}
\item[(a)] $\lambda(z) \, |dz|$ is complete.
\item[(b)] Every branch of the inverse of the developing map $f : G \to \D$ can be
  analytically continued along any path $\gamma \subset \D$.
\item[(c)] $\lambda(z) \, |dz|$ is the hyperbolic metric of $G$.
\end{itemize}
\end{theorem}

{\bf Proof.}

(a) $\Longrightarrow$ (b):
Let $\lambda(z) \, |dz|$ be a complete conformal metric 
 of constant curvature $-1$ in $G$. By Liouville's Theorem  we have
$$ \lambda(z) \, |dz|=\frac{2\, |f'(z)|}{1-|f(z)|^2} \, |dz|$$
for some locally univalent (multivalued) function $f : G \to \D$.
Without loss of generality we may assume $0 \in G$ and $f(0)=0$. 
We claim that the branch $f^{-1}$ of the inverse function of $f$ 
with $f^{-1}(0)=0$ can be continued analytically along every path 
$\gamma : [a,b] \to \D$ with $\gamma (a)=0$.

\medskip

If this were false, then there exists a number $\tau \in (a,b)$ such
that $f^{-1}$ can be continued analytically along
$\gamma : [a,\tau) \to \D$, but not further. Let $\tau_n \in [a,\tau)$
with $\tau_n \to \tau$ and let $\gamma_n:=\gamma|_{[a,\tau_n]}$.
Note that 
$$ \int \limits_{\gamma_n} \lambda_{\D}(w) \, |dw| \le
\int \limits_{\gamma} \lambda_{\D}(w) \, |dw| \le C$$
for some constant $C>0$. Let $w_n=\gamma(\tau_n)$ and
$z_n=f^{-1}(w_n)$. Then we get
$$ \text{d}_{\lambda}(z_n,0) \le \int \limits_{f^{-1} \circ
  \gamma_n} \lambda(z) \, |dz|=\int \limits_{\gamma_n} \lambda_{\D}(w)
\, |dw| \le C\, .
$$
As $\lambda(z) \, |dz|$ is
complete for $G$, the sequence $(z_n)$
 stays in some compact subset of $G$, so some subsequence converges to 
a point $z_{\infty} \in G$. We may thus assume $z_n \to z_{\infty} \in G$. Hence
$\gamma (\tau) \leftarrow w_n=f(z_n) \rightarrow f(z_{\infty})$,
i.e.~$f(z_{\infty})=\gamma(\tau)$. Since $f$ is locally univalent, 
it maps a neighborhood of
$z_{\infty}$
univalently onto a disk $K_r(\gamma(\tau))$, that is 
$f^{-1}$ can be continued analytically to this disk.
This, however, contradicts our hypothesis.

\medskip

(b) $\Longrightarrow$ (c):
Note that $g:=f^{-1}: \D \to G$ is  holomorphic and single--valued by the
Monodromy Theorem. Hence  the Fundamental Theorem gives us 
$$ (g^{\, *}\lambda_G) (z) \le \lambda_{\D}(z)\,,  \quad z \in \D\,.$$
But, by construction, 
$$ \lambda_{\D} (z)= (g^{ *}\lambda)(z)\,,  \quad z \in \D\,.$$
Thus $\lambda_G \le \lambda$  in $G$ and therefore
  $\lambda\equiv \lambda_G$ in $G$.

\medskip

(c) $\Longrightarrow$ (a): 
This is Theorem \ref{thm:hyp20}. \hfill{$\blacksquare$}



\begin{theorem}[Uniformization Theorem]
Let $G$ be a domain with at least two boundary points. Then there exists a locally univalent, surjective and  analytic function $\pi: \D \to G$ such that every branch of $\pi^{-1}$ can be analytically continued along any path in $G$ ($\pi$ is called a universal covering of G). If $\tau$ is another universal covering of G, then $ \tau= \pi \circ T$ for some conformal self--map $T$ of $\D$. 
\end{theorem}

{\bf Proof.} $G$ carries an SK--metric, see Theorem \ref{thm:hyp}, and
thus possesses a hyperbolic metric $\lambda_G(w) \, |dw|$ by Theorem \ref{thm:perron}, which is
complete by Theorem \ref{thm:hyp20}. Theorem \ref{thm:liouvillecomplete}
implies that the every branch $\pi$ of the inverse 
of the developing map $f : G \to \D$ can be analytically continued along any
path in $\D$. By the Monodromy Theorem $\pi$ is an analytic function from
$\D$ into $G$. It follows from the construction that $\pi$ is locally
univalent. Since $f$ can also be analytically continued
along any path in $G$ (see Remark \ref{rem:21}) it follows that $\pi$ is onto.
 If $\tau$ is another such function, then  it is easy to see that $\tau=\pi \circ T$ for some
conformal self--map $T$ of $\D$. \hfill{$\blacksquare$}

\begin{remark1}
Let $\lambda_G(z)\, |dz| $ be the hyperbolic metric on $G$ and $\pi : \D \to
G$ a universal covering. 
 Then $( \pi^* \lambda_G)(z)\, |dz|= \lambda_{\D}(z) \, |dz|$ by construction. 
\end{remark1}

\begin{corollary}[Riemann Mapping Theorem]
Let $G \subsetneq \C$ be a simply connected domain. Then there exists a conformal map
$f$ from $G$ onto $\D$, which is uniquely determined up to postcomposition
with a conformal self--map of $\D$.
\end{corollary}

{\bf Proof.} If $G$ is simply connected, then the developing map $f$ is a
single--valued analytic function and therefore a conformal map from $G$ onto $\D$.
 \hfill{$\blacksquare$}

\subsection{The twice--punctured plane}

\medskip

In case of the twice--punctured plane it is possible to derive explicit formulas for the hyperbolic metric and its developing map.

\begin{theorem}[Agard's formula] \label{thm:agard}
Let $$ K(z) :=\frac{2}{\pi} \int \limits_0^1 \frac{dt}{\sqrt{(1-t^2) \, (1-z t^2)}} \, .$$
Then 
$$ \lambda_{\C''}(z)\, |dz|=\frac{|dz|}{\pi \, |z| \, |1-z| \, \Re \left[ K(z) K(1-\overline{z})\right]}  $$
and  the developing map of $\lambda_{\C''}(z) \, |dz|$
is given by
 $$ z \mapsto \frac{K(1-z)-K(z)}{K(1-z)+K(z)}\, . $$
\end{theorem}

The proof of Theorem \ref{thm:agard} relies on Liouville's Theorem and the
following lemma.

\begin{lemma} \label{lem:aga1}
Let $\lambda(z) \, |dz|=\lambda_{\C''}(z) \, |dz|$ be the hyperbolic metric of $\C''$. Then
$$ S_{\lambda}(z)=\frac{1}{2} \left[ \frac{1}{z^2}+\frac{1}{(z-1)^2}+\frac{1}{z \, (1-z)} \right] \, .$$
\end{lemma}

{\bf Proof.}
Theorem \ref{thm:sch} shows that $S_{\lambda}$ is analytic in $\C''$ with poles
of order $2$ at $z=0$ and $z=1$, so 
$$ S_{\lambda}(z)=\frac{1}{2 \, z^2}+\frac{c_1}{z}+\frac{1}{2 \, (z-1)^2}
+\frac{c_2}{z-1}+r(z)=\frac{1}{2 \, z^2}+\frac{1}{2 \, (z-1)^2}
+\frac{(c_1+c_2) z-c_1}{z (z-1)}+r(z) \, ,$$
where $c_1, c_2 \in \C$ and $r$ is analytic on $\C$.
We analyze $S_{\lambda}(z)$ at $z=\infty$.
Let $f(z)=1/z$ and observe that $f^*\lambda=\lambda$, so
$S_{\lambda}(z)=S_{f^*\lambda}(z)=S_f+S_{\lambda}(f(z)) \, f'(z)^2=S_{\lambda}(1/z)/z^4$ in view of Exercise \ref{sec:constcurv}.\ref{exe:schwarztransform}.
Hence 
$$ \lim \limits_{z \to \infty} z^2 S_{\lambda}(z)=\lim \limits_{z \to \infty} \frac{1}{z^2} S_{\lambda}(1/z)=\lim \limits_{z \to 0} z^2 S_{\lambda}(z)=\frac{1}{2} \, .$$
 This forces $r \equiv 0$ 
 and $c_1=-c_2=1/2$. 
\hfill{$\blacksquare$}

\medskip

{\bf Proof of Theorem \ref{thm:agard}.}
We consider the simply connected domain $G=\C\backslash ((-\infty,0] \cup [1,+\infty))$. Liouville's Theorem gives us an analytic function $f : G \to \D$ such that
$$ \frac{2 \, |f'(z)|}{1-|f(z)|^2}=\lambda_{\C''}(z) \, , \qquad z \in G \, . $$
By Example \ref{l1} and Lemma \ref{lem:aga1}, we know that $f$ is a solution to
\begin{equation} \label{eq:s34}
S_f(z)=\frac{1}{2} \left[ \frac{1}{z^2}+\frac{1}{(z-1)^2}+\frac{1}{z \, (1-z)} \right]\,. 
\end{equation}
We consider the hypergeometric differential equation
\begin{equation} \label{eq:hypgeometric}
 z (1-z) \, u''+\left( 1-2 \,z\right) \, u'-\frac{1}{4} \, u=0 \, ,
\end{equation}
since  every solution of (\ref{eq:s34}) can be written as a fractional linear 
transformation of two linearly independent  solutions $u_1$ and $u_2$ of
(\ref{eq:hypgeometric}), see \cite[p.~203 ff.]{Neh52}. It is well--known and easy to prove  (see Nehari \cite[p.206 ff.]{Neh52}) that 
 $u_1(z)=K(z)$, $u_2(z)=K(1-z)$ are two linearly independent  solutions 
of (\ref{eq:hypgeometric}) in $G$. Hence 
$$ f(z)=\frac{a \, u_2(z)+b \, u_1(z)}{c \, u_2(z)+d \, u_1(z)}$$
for $a,b,c,d \in \C$ with $ad-bc\not=0$ and thus
$$  \lambda_{\C''}(z)=2 |ad-bc| \, \frac{|u_2'(z) u_1(z)-u_2(z) u_1'(z)|}{|c \,u_2(z)+
d \, u_1(z)|^2-|a \, u_2(z)+b \, u_1(z)|^2} \, , \qquad z \in G \, . $$
Now, consider the Wronskian $w:=u_2'u_1-u_2 u_1'$  of $u_1$ and $u_2$, which 
is a solution  to 
$$ w'=-\frac{1-2 z}{z \, (1-z)} \, w \, . $$
As  $u_1$ is holomorphic in a neighborhood of $z=0$ with $u_1(0)=1$ and 
$$u_2(z)=-\frac{1}{\pi} \log z + h_1(z) \, \log z+ h_2(z)\, ,$$
where $h_1$ and $h_2$ are analytic in a neighborhood of $z=0$ with
$h_1(0)=0$, cf.~\cite[15.5.16 and 15.5.17]{AB}), we obtain
 
$$ u_2'u_1-u_2 u_1'=- \frac{1}{\pi}\,  \frac{1}{z \, (1-z)} \, . $$
Consequently, if we set $\beta:= |ad-bc|$ and $\gamma:=\overline{c} d-b\overline{a}$, we get
\begin{equation*} \label{eq:aga2}
 \lambda_{\C''}(z)= \frac{2\, \beta}{ \pi } \frac{1}{ |z| \, |1-z|} \frac{1}{\left( |c|^2-|a|^2 
\right) |u_2(z)|^2+\left(
|d|^2-|b|^2 \right)  |u_1(z)|^2+2 \Re \big( \gamma \,  u_1(z) \overline{u_2(z)} \big) } \, , 
\end{equation*}
for every $z \in G$.  Now, if we let $z \to 0$, then by the  asymptotics of
$\lambda_{\C''}(z)$ at $z=0$ described in Corollary \ref{cor:is} we see that
$|a|=|c|$. A similar analysis for $z \to 1$ gives $|b|=|d|$. Hence 
\begin{equation} \label{eq:89}
 \lambda_{\C''}(z)=\frac{\beta}{\pi }\,  \frac{1}{|z| \, |1-z| \,  \Re \left(
     \gamma \, u_1(z) \,  \overline{u_2(z)}  \right) } \, .
\end{equation}


Applying  Corollary \ref{cor:is} once more shows $\beta=
\Re(\gamma)$. In addition, the continuity of $\lambda_{\C''}(z)$  in $\C''$
implies that $\gamma$ is real.  We thus arrive at
$$ \lambda_{\C''}(z)\, |dz|=\frac{|dz|}{\pi \, |z| \, |1-z| \, \Re \left[ K(z) K(1-\overline{z})\right]} \, . $$

Now, we are left to show that some developing map $F : \C'' \to \D$ of $\lambda_{\C''}(z) \, |dz|$ has the asserted form.
Here we will make essential use of the fact that with $f$ also $T \circ f$ is 
a developing map for $\lambda_{\C''}(z) \, |dz|$, when $T$ is a conformal
self--map of $\D$. We have shown above that
$$f(z)= \frac{a\, u_2(z) +b \, u_1(z) }{c \, u_2(z) + d u_1(z)}$$
with $|a|=|c|$, $|b|=|d|$ and $\overline{c} d-b\overline{a}=|ad -bc| >0$ is a
developing map for $\lambda_{\C''}(z) \, |dz|$. Therefore,
$$\tilde{F}(z)=\frac{u_2(z) + B\, u_1(z)}{ u_2(z) + D\, u_1(z)}$$
 is also a developing map for $\lambda_{\C''}(z) \, |dz|$ where $B=b/a$ and $D=d/c$. Furthermore, $|B|=|D|$ and
$D-B>0$. This implies $\Re(B)<0<\Re(D)$. Thus, $z_0:=(1+B)/(1+D) \in \D$ 
and if 
$$T(z)=  \frac{1-\overline{z_0}}{1- z_0 }\, \frac{z -z_0}{1-\overline{z_0}
  z}\, , $$
then
$$F(z) =(T \circ \tilde{F})(z) = \frac{u_2(z) -u_1(z)}{u_2(z)+u_1(z)}= \frac{K(1-z)-K(z)}{K(1-z)+K(z)} \,  $$
is the desired developing map of $\lambda_{\C''}(z) \, |dz|$.
\hfill{$\blacksquare$}

\begin{corollary} \label{cor:val}
The density of the hyperbolic metric $\lambda_{\C''}(z) \, |dz|$ at $z=-1$
has the value
$$\lambda_{\C''}(-1)=\frac{\Gamma \left( \frac{3}{4} \right)^4}{\pi^2} \approx \,0.22847329 . $$
\end{corollary}

{\bf Proof.}
Since $K(1/2)=\sqrt{\pi}/\Gamma(3/4)^2$ (see \cite{AB}), 
Exercise \ref{sec:hyperbolic}.\ref{exe:wert1/2} and
Theorem \ref{thm:agard} give
$$ \lambda_{\C''}(-1)=\frac{\lambda_{\C''}(1/2)}{4} =\frac{1}{\pi \, K(1/2)^2} =
\frac{\Gamma \left( \frac{3}{4} \right)^4}{\pi^2}
\, . $$
\hfill{$\blacksquare$}\\
\begin{figure}[h]
\centerline{\includegraphics[width=12cm]{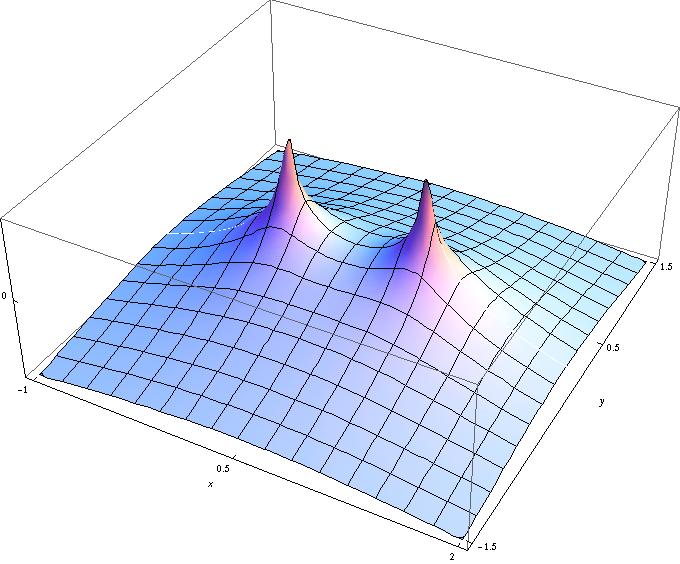}}
\caption{A plot of $z \mapsto \log \lambda_{\C''}(z)$}
\end{figure}

\section*{Exercises for Section 5}

\begin{enumerate}
\item \label{exe:5.1} Let $\lambda(z) \, |dz|$ be a regular conformal metric with constant curvature $+1$ on a simply connected domain $G$. Show that
$$ \lambda(z)=\frac{2 \, |f'(z)|}{1+|f(z)|^2}$$
for some locally univalent meromorphic function $f$ on $G$.
\item \label{exe:schwarztransform}
Let $\lambda(w) \, |dw|$ be a regular conformal pseudo--metric on $G$ and
$f : D \to G$ an analytic map. Show that
$$ S_{f^*\lambda}(z)=S_f(z)+S_{\lambda}(f(z)) \, f'(z)^2\,. $$
\end{enumerate}

\section*{Notes}
Theorem \ref{thm:liouville} was apparently first stated by Liouville 
\cite{Lio1853}, but his proof is certainly not complete by today's standards. The proof given here is an adaption of the method of Nitsche 
\cite{Nit57}, whose presentation contains some gaps as well.
Liouville's theorem is also proved in Bieberbach \cite{Bie16} and a very elegant geometric proof has been given by D.~Minda (unpublished).
Theorem \ref{thm:sch} is Theorem 1 (iii) in Minda \cite{Min97}, but can also be
deduced from Nitsche's result in \cite{Nit57}. The methods of Minda and 
Nitsche are different from the short proof given above.


\section{Appendix}

In this appendix we construct conformal metrics  
with constant curvature $-1$ on disks with prescribed (continuous) boundary values
(Theorem \ref{thm:existence}). It clearly
suffices to consider the case of the unit disk $\D$. In this situation,
Theorem \ref{thm:existence} is equivalent to the following result.

\begin{theorem} \label{solgen}
Let  $\psi : \partial \D \to
\R$ be a continuous function.
Then there exists a uniquely determined function $u \in C(\overline{\D} ) \cap C^2(\D )$ such that \label{ref:diff} \label{def:C}
\begin{equation} \label{eq:pdesing}
 \begin{array}{rccc}
\Delta u &=&  e^{2 u} & \, \, \text{in }\, \,   \D, \\[2mm]
u &=& \psi & \, \, \text{on } \, \,  \partial \D.
\end{array}
\end{equation}
\end{theorem}

\medskip

We shall use the following basic facts about the Poisson equation 
$\Delta u=q$ in $\D$.

\begin{theorem} \label{thm:Gilbarg}
 Let $q : \D \to \R$ be a bounded and continuous function. Then 
\begin{equation*}
w(z):=  -\frac{1}{2 \pi} \iint \limits_{\D} g(z, \zeta) \, q(\zeta) \, d
m_{\zeta}
\end{equation*}
where $g(z, \zeta)=\log \left(\left|1 - \overline{\zeta} z\right| / \left|z- \zeta
\right|\right) $ is Green's function for $\D$, 
 belongs to $C(\overline{\D}) \cap C^1(\D)$. \label{ref:Diff}
If, in addition,  $q$ is continuously differentiable in $\D$  then $w \in
C^{2}(\D)$ and solves the boundary value problem
\begin{equation*}
 \begin{array}{rccc}
\Delta w &=&  q  &\, \, \text{in } \, \,  \D, \\[2mm]
w &=& 0 & \, \, \text{on } \, \,  \partial \D.
\end{array}
\end{equation*}
Moreover, if $q \in C^k(\D)$ for some $k \ge 1$, then $w \in C^{k+1}(\D)$.
\end{theorem}

\begin{remark}\label{rem:representation}
Theorem \ref{thm:Gilbarg} shows that every
$C^2$ solution to (\ref{eq:pdesing}) has the form
\begin{equation} \label{eq:form}
u(z)=h(z)-\frac{1}{2 \pi}\iint \limits_{\D} g(z,\zeta) \,
e^{2 u(\zeta)} \, dm_{\zeta}\,, \qquad z \in \D\, , 
\end{equation}
where $h$ is continuous on $\overline{\D}$, harmonic in $\D$  and coincides with 
$\psi$ on $\partial \D$. 
If, conversely,  a bounded and continuous function $u$ on $\D$ has the  form
(\ref{eq:form}) with a harmonic function $h$  in $\D$ which is continuous on $\overline{\D}$
and has boundary values $\psi$, then $u \in
C^2(\D)$ and solves  (\ref{eq:pdesing}).
In particular, every $C^2$ solution to $\Delta u=e^{2u}$ is of class
$C^{\infty}$.
\end{remark}

{\bf Proof of Theorem \ref{thm:Gilbarg}.}
For the differentiability properties of $w$ we refer to \cite[p.~50
ff]{GT77}. So, it remains to show that
$$ \lim \limits_{z \to \xi} \iint \limits_{\D} g(z,\zeta) \, q(\zeta) \,
dm_{\zeta}=0\, .
$$

Since $\zeta \mapsto g(z, \zeta)$ has no integrable majorant independent of $z$, we
cannot make profit  of the fact that $z \mapsto g(z, \zeta)$ vanishes continuously on the
boundary of $\D$. We need a more refined  argument. For this we
 pick $\tau \in  \partial \D$ and choose $\eps >0$. Then we can find $0<
\varrho < 1/4$ such that for any $z\in \D \cap K_{\varrho/2}(\tau)$ 
\begin{eqnarray*}
\displaystyle \iint \limits_{\D \cap K_{\varrho} (\tau)} \log \left| \frac{1-
      \overline{\zeta}z}{z- \zeta}\right| \, dm_{\zeta} & \le & (\log 2)\,  \pi \,
  \varrho^2 + \iint \limits_{ K_{\varrho} (\tau)} \log \frac{1}{|z-
    \zeta|}\, dm_{\zeta}    \\[1mm]
\displaystyle & \le & (\log 2)\,  \pi \,\varrho^2  + \iint \limits_{ K_{2
      \varrho} (z)} \log \frac{1}{|z- \zeta|} \, dm_{\zeta}   \le (\log 2)\,  \pi
\,\varrho^2  +4\, \pi \, \varrho \le \eps \, .
\end{eqnarray*}

If $z \in \D \cap K_{\varrho/2}(\tau)$ and $\zeta \in \D \backslash
K_{\varrho}( \tau)$ then
\begin{equation*}
\displaystyle{\left| \frac{1-  \overline{\zeta}z}{z- \zeta}    \right| \le 
 \frac{|1-  \overline{\zeta} \tau| + | z- \tau|}{|\tau- \zeta|- |z -\tau|}
 = \frac{1 + \frac{| z- \tau|}{|\tau- \zeta|} }{1- \frac{|z
    -\tau|}{|\tau- \zeta|} } \le \frac{1 + \frac{| z- \tau|}{\varrho}
}{1- \frac{|z
    -\tau|}{\varrho} } }\, .
\end{equation*}
This implies that
$$ \iint \limits_{\D \backslash K_{\varrho} (\tau)} \log \left| \frac{1-
      \overline{\zeta}z}{z- \zeta}\right| \, dm_{\zeta} \le  \pi \, \log
  \left( \frac{\varrho + |z - \tau|}{\varrho -|z - \tau|}
   \right)  \le \eps $$
if $z \in \D \cap K_{\varrho/2}(\tau)$ is sufficiently close to $\tau$. 
The result follows.~\hfill{$\blacksquare$}

\medskip

The next lemma provides an important equicontinuity property of Green's function.

\begin{lemma}\label{lem:equicontinuous}
Let $0< \varrho<1 $. Then for every $z_1,z_2 \in
\overline{\D}_{1-\varrho}$ with $|z_1-z_2|< \varrho/4$,
\begin{equation}\label{uniform}
 \iint \limits_{\D}\,  |g(z_1, \zeta)-g(z_2, \zeta)| \, dm_{\zeta} \le   
 \, \pi \, ( 5/\varrho+4 \, \varrho) \, |z_1-z_2|\, .
\end{equation}
\end{lemma}

{\bf Proof.}
Let $z_1,z_2 \in \overline{\D}_{1- \varrho}$ such that $|z_1-z_2| < \varrho/4
  $. We first observe that
\begin{equation*}
|g(z_2, \zeta)-g(z_1,\zeta)| 
 \le \displaystyle \left| \log \left| \frac{z_1- \zeta}{z_2 - \zeta}\right|  \right|   + \left| \log  \left| \frac{1-\overline{\zeta}z_2}{1- \overline{\zeta}z_1} \right|  \right| \displaystyle \le \left| \log\left| \frac{z_1- \zeta}{z_2 - \zeta}\right|  \right| + \log
 \left(1+\frac{|z_1-z_2|}{\varrho} \right) \, .
\end{equation*}
Thus  we have
$$\iint \limits_{\D} |g(z_1, \zeta)-g(z_2,\zeta)| \, dm_{\zeta} \le \pi  \,
\frac{|z_1-z_2|}{\varrho} + \iint \limits_{\D} \left| \log\left|
     \frac{z_1- \zeta}{z_2 - \zeta}\right|  \right| \, dm_{\zeta}\, .$$

Now let $M=(z_1+z_2)/2$ be the midpoint of the line segment joining 
$z_1$ and $z_2$ and let
$K=K_{ \varrho /2}(M)$. Then we obtain
\begin{equation*}
\begin{array}{l}
\displaystyle{\iint \limits_{\D} \left| \log\left|
     \frac{z_1- \zeta}{z_2 - \zeta}\right|  \right| \, dm_{\zeta}}  = \displaystyle{  \iint
 \limits_{K } \left| \log\left|
     \frac{z_1- \zeta}{z_2 - \zeta}\right|  \right| \, dm_{\zeta} + \iint
 \limits_{\D \backslash K } \left| \log\left|
     \frac{z_1- \zeta}{z_2 - \zeta}\right|  \right| \, dm_{\zeta} }\\[8mm]
 \quad \le  \displaystyle{  \iint
 \limits_{K }  \left[ \log \left( 1 + \frac{|z_2-z_1|}{|z_1 - \zeta|}   \right) +
 \log \left( 1 + \frac{|z_2-z_1|}{|z_2 - \zeta|}   \right) \right] dm_{\zeta}  + \iint
 \limits_{\D \backslash K } \log \left( 1+ \frac{4\, |z_1 -z_2|}{\varrho}
 \right) dm_{\zeta}}\\[8mm]
\quad \le \displaystyle{ \iint \limits_{K_{\varrho }(z_1)}  \log \left( 1 +
  \frac{|z_2-z_1|}{|z_1 - \zeta|}   \right)\, dm_{\zeta}+ \iint\limits_{K_{\varrho}(z_2)}
\log \left( 1 + \frac{|z_2-z_1|}{|z_2 - \zeta|}   \right) \, dm_{\zeta}+ 
\frac{4}{\varrho} \, \pi |z_1-z_2| }\\[8mm]
\quad \displaystyle{\le 4\, \pi \, \varrho\, |z_1-z_2| \, + \,  \frac{4}{\varrho} \,
  \pi |z_1-z_2| }\, . \hfill \blacksquare
\end{array}
\end{equation*}

\medskip

We are now in a position to give the

\medskip

{\bf Proof of Theorem \ref{solgen}.}

{\bf Uniqueness:}
Assume $u$ and $v$ are two solutions of the boundary value
problem (\ref{eq:pdesing}). Then $s(z):= \max\{0, u(z)-v(z)\}$ is a
non--negative subharmonic function in $\D$ and $s(z)=0$ on $|z|=1$. By the
maximum principle for subharmonic functions it follows that $u(z) \le v(z)$ for every
$z \in \D$. Switching the r$\hat{\mbox{o}}$les of $u$ and $v$ gives $u(z)=v(z)$ for $z \in
\D$.

\medskip

{\bf Existence:}
Suppose for a moment $u \in C(\overline{\D})
\cap C^2(\D)$ is a solution to (\ref{eq:pdesing}). 
Then Theorem \ref{thm:Gilbarg} 
implies that 
$$u(z)=h(z)-\frac{1}{2 \pi}\iint \limits_{\D} g(z,\zeta) \,
e^{2 u(\zeta)} \, dm_{\zeta}\,, \qquad z \in \D\, , $$
where $h$ is continuous on $\overline{\D}$, harmonic in $\D$  and coincides with 
$\psi$ on $\partial \D$. 

\medskip

This suggests to introduce the operator
$$ T[u](z):= h(z)-\frac{1}{2 \pi}\iint \limits_{\D} g(z,\zeta)\, 
e^{2 u(\zeta)} \, dm_{\zeta},$$
and to apply Schauder's fixed point theorem.

\medskip
To set the stage for Schauder's theorem,
let $X$ be the Fr\'echet space of all real--valued continuous functions
in $\D$ equipped with the (metriziable) compact--open topology, and let
$$M:=\{ u\in X \, : \, m \le u(z) \le h(z) \, \, \text{for all}\, z \in \D
\}\,, \text{ where }  \,  m:=\inf_{z \in \D} T[h](z)\, .$$
Note that $m>-\infty$.
In order to be able to apply  Schauder's fixed
point\footnotemark \footnotetext{compare \cite[p. 90]{Dei85}}
theorem we need to check the following properties:

\medskip

 The set $M$ is closed and convex (in $X$). The
  operator $T : M \to X$ is continuous, maps $M$ into $M$, and
$T(M)$ is precompact.

\medskip

Clearly,  $M$ is closed and convex.
Next, we will prove $T[M]$ is precompact, by showing $T[M]$ is a locally
 equicontinuous family  and the set $\{Tu(z) \, : \, u \in M  , \, z \in \D\}$
is bounded in the reals. For that
 pick $0< \varrho <1$,  let $B_{1-\varrho}:= \{ z \in \C: |z| \le 1- \varrho \}$
 and  fix $\eps>0$. 
Since $h$ is continuous on $\overline{\D}$ there exists a constant
$\delta'>0$ such that $|h(z_1)-h(z_2)|<\eps/2$ for all $z_1, z_2 \in
\D$ with $|z_1-z_2|<\delta'$. We now define
$$ \delta:=\min \left\{ \delta', \frac{  \eps}{ C \left( 5/\varrho +4
      \varrho \right)}, \frac{\varrho}{4} \right\}  \, ,$$
where $C:= \sup_{\zeta \in \D}   e^{2 h(\zeta)}$.
 Then by Lemma \ref{lem:equicontinuous} we have  for all
$z_1,z_2 \in B_{1-\varrho}$ with $|z_1-z_2|<\delta$ and for every $u \in M$: 
\begin{eqnarray*} |T[u](z_2)-T[u](z_1)| & \le &
|h(z_2)-h(z_1)|+\frac{1}{2 \pi} \iint \limits_{\D} \left|
 g(z_2,\zeta)-g(z_1,\zeta)\right| \, e^{2 u(\zeta)} \,
 dm_{\zeta}\\
& \le & |h(z_2)-h(z_1)|+ \frac{1}{2 \pi} \, C \iint \limits_{\D} \left|
 g(z_2,\zeta)-g(z_1,\zeta)\right|  \,
 dm_{\zeta}\\
& \le & \frac{\eps}{2}+ \frac{1}{2 \pi} \, C\,  \pi (5/\varrho + 4\, \varrho) 
 \cdot \delta \le \eps.
\end{eqnarray*}
Thus $T[M]$ is a locally equicontinuous set of 
functions on $\D$. Moreover, for all $u \in M$ and all $z \in \overline{\D}$ 
\begin{equation} \label{oli}
T[h](z) \le T[u](z) \le h(z) \, , 
\end{equation}
which implies
$$ \min \limits_{\zeta \in \overline{\D}} T[h](\zeta) \le T[u](z) \le
\max \limits_{\zeta \in \overline{\D}} |h(\zeta)| \quad \text{for every} \, u
\in M\, \, \text{and for all} \, z \in \overline{\D}\,
.$$
This shows $\{Tu(z) \, : \, u \in M, \, z \in \D\}$
 is bounded, so $T[M]$ is a precompact subset of $X$.
Note, estimate (\ref{oli}) also gives $T[M] \subseteq M$. 

\medskip 
It  remains to prove that
 $T : M \to M$ is continuous. 
Let $(u_k)$ be a sequence of functions in $M$ which converges
locally uniformly in $\D$ to $u \in M$. We have to show that
the sequence $(T[u_k])$ converges locally uniformly in $\D$ to
$T[u]$.

\medskip

To do this, choose  $0<\varrho<1$  and fix $\eps>0$. We shall
deduce that $|T[u_k](z)-T[u](z)|<\eps$ for all $z \in B_{1-\varrho}$ and all
$k \ge \tilde{k}$ for some $\tilde{k}$ independent of $z$.

\medskip

For notational simplicity we set

$$ C_1:= \frac{1}{ \pi} \sup  \limits_{\zeta \in \D} e^{2
  h(\zeta)}\, , \qquad C_2 := \sup_{z \in \D}\iint \limits_{\D} g(z, \zeta) \, dm_{\zeta}\, , \qquad
C_3 := 2 \sup \limits_{\zeta \in \D} |h(\zeta)|, 
$$

Now we choose $0<r <\varrho/2$ such that
\begin{equation}\label{eq:ex1}
C_1\, C_3 \, \log \left(\frac{2}{r} \right) \, \pi\, (2r -r^2) < \frac{\eps}{2}\, .\end{equation}
Further, we can find  an index $\tilde{k} \in \N$ such that
\begin{equation}\label{eq:ex2}
 \sup_{z \in B} \left|u_k(z)-u(z) \right| <\frac{\eps}{2 C_1 C_2}, \qquad
k \ge \tilde{k}\, ,
\end{equation}
where $B:=  B_{1-r}$. We now obtain by equations (\ref{eq:ex1}) and (\ref{eq:ex2}) for $z \in B_{1-\varrho}$  

\begin{equation*}
\begin{array}{l}
\displaystyle\left| T[u_k](z)-T[u](z) \right| = \left| \frac{1}{2\pi} \iint \limits_{\D}g(z,
  \zeta) \,  (e^{2u_k(\zeta)} -e^{2u(\zeta)}) dm_{\zeta}\right|\\[6mm]
 \qquad  \displaystyle\le   \frac{1}{ \pi} \,  \sup \limits_{\zeta \in \overline{\D}} e^{2
  h(\zeta)}\cdot \iint \limits_{\D} g(z,\zeta) \left|
  u(\zeta)-u_k(\zeta) \right| \, dm_{\zeta} \\[2mm] 
 \qquad   \displaystyle \le   C_1  \sup_{z \in B} \left|u_k(z)-u(z) \right| \iint \limits_{\D} g(z,
\zeta) \, dm_{\zeta} \,  + \, C_1 \,  2 \, \sup_{\zeta \in \D} |h(\zeta)| \iint
\limits_{\D \backslash B} g(z, \zeta)  \, dm_{\zeta} \\
\qquad  \displaystyle \le  \frac{\eps}{2} + C_1 \, C_3   \iint \limits_{\D \backslash B}
  \log\left(\frac{2}{r} \right)  \,
dm_{\zeta}< \eps   \, .
\end{array}
\end{equation*}
Thus $T:M \to M$ is continuous.

\medskip

Now Schauder's fixed point theorem  gives us  a fixed
point $u \in M$ of $T$, i.e.,
\begin{equation*} \label{eq:fixedpoint}
 u(z)=T[u](z)=h(z)-\frac{1}{2 \pi} \iint \limits_{\D} g(z,\zeta) \,  e^{2
   u(\zeta)} \, dm_{\zeta}.
\end{equation*}
 We claim $u$ belongs to
$C(\overline{\D})\cap C^2(\D)$ and is a solution of (\ref{eq:pdesing}).

\medskip

Indeed, since $\zeta \mapsto  e^{2 u(\zeta)}$ is
bounded and continuous in $\D$, the function $u$ belongs to $
C(\overline{\D})\cap C^1(\D)$ by Theorem \ref{thm:Gilbarg}. 
This implies that the function $\zeta \mapsto  e^{2 u(\zeta)}$ belongs to
$C^1(\D)$ and
applying Theorem \ref{thm:Gilbarg} again proves that $u \in C(\overline{\D}) \cap
C^2(\D)$ solves (\ref{eq:pdesing}).
\hfill{$\blacksquare$}

\end{document}